\documentclass[11pt]{amsart}
\usepackage{amsthm}
\usepackage{amssymb}
\usepackage{amsmath}

\newcommand     {\comment}[1]   {}
\newcommand{\mute}[2] {}
\newcommand     {\printname}[1] {}
\newcommand{\labell}[1] {\label{#1}\printname{#1}}

\numberwithin{equation}{section}
\newtheorem {Theorem}			{Theorem}

\newtheorem {Lemma}[equation]     	{Lemma}

\newtheorem {Proposition} [equation]	{Proposition}

\theoremstyle{definition}
\newtheorem{Definition}[equation]{Definition}

\theoremstyle{remark}
\newtheorem{Remark}[equation]{Remark}
\newtheorem*{remark}{Remark}
\newtheorem{Example}[equation]{Example}

\def	\Z	{{\mathbb Z}}
\def	\R	{{\mathbb R}}
\def	\C	{{\mathbb C}}

\def	\W	{{\mathcal W}}

\def	\t	{{\mathfrak t}}

\def	\h	{{\mathfrak h}}
\def	\k	{{\mathfrak k}}

\newcommand\fh {{\mathfrak h}}
\newcommand\fk {{\mathfrak k}}
\newcommand\ft {{\mathfrak t}}

\def	\bfH	{\mathbf{H}}
\def	\bfG	{\mathbf{G}}

\newcommand\fU {{{\mathfrak U}}}
\newcommand\fV {{{\mathfrak V}}}

\def	\fh	{{\mathfrak h}}

\def	\cQ		{{\mathcal Q}}
\def	\hcQ		{{\hat{\mathcal Q}}}
\def	\cRQ		{{\mathcal R}{\mathcal Q}}
\def	\cP		{{\mathcal P}}
\def	\hcP		{{\hat{\mathcal P}}}
\def	\cRP		{{\mathcal R}{\mathcal P}}
\def	\cE		{{\mathcal E}}

\def    \hcH {{\hat{\mathcal H}}}

\def	\PhiT		{\Phi \mbox{-} T}

\def	\psibar  {{\ol{\psi}}}

\def	\Phibar  {{\ol{\Phi}}}
\def	\Rbar	{{\ol{R}_Y}}

\newcommand\cA{{\mathcal A}}
\newcommand\cB{{\mathcal B}}
\newcommand\cF{{\mathcal F}}

\def \calT {{\mathcal T}}

\newcommand\cO{{\mathcal O}}

\def	\vZ	{\text{\v{Z}}}
\def	\vH	{\operatorname{H}}

\def    \ob    {{\operatorname{ob}\,}}
\def    \hom    {{\operatorname{hom}}}
\def    \exc    {{\operatorname{exc}}}

\def	\GL	{{\operatorname{Gl}}}

\def	\Stab	{\operatorname{Stab}}

\def	\image	{\operatorname{image}}

\def	\GL	{\operatorname{GL}}

\def	\SU	{\operatorname{SU}}
\def	\SO	{\operatorname{SO}}

\def	\max	{{\operatorname{max}}}
\def	\min	{{\operatorname{min}}}

\def \cutoff {{\text{cutoff}}}

\def	\ss	{\scriptstyle}
\def	\bs	{\boldsymbol}
\def	\inv	{^{-1}}
\def	\to	{\longrightarrow}
\def	\ssminus	{\smallsetminus}

\def	\half	{\frac{1}{2}}
\def	\Cc	{{\C^\times}}
\def	\DH	{Duistermaat-Heckman\ }

\def	\id	{\text{id}}

\newcommand     {\ol}[1]   {\overline{#1}}

\begin{document}

\title[Hamiltonian torus actions with two dimensional quotients]
{Complete invariants for Hamiltonian torus actions with two dimensional 
quotients}

\author[Yael Karshon]{Yael Karshon}
\address{Inst.\ of Mathematics,
The Hebrew University of Jerusalem, Israel,
and: Dept.\ of Mathematics, the University of Toronto,
Toronto, Ontario M5S 3G3, Canada}
\email{karshon@math.toronto.edu}

\author[Susan Tolman]{Susan Tolman}
\address{Dept.\ of Mathematics, The University of Illinois
at Urbana-Champaign, 1409 W.~Green Street, Urbana, IL 61801, U.S.A}
\email{stolman@math.uiuc.edu}

\thanks{
S. Tolman was partially supported by a Sloan fellowship and by
NSF grants DMS-980305 and DMS 02-04448.
Both authors were partially supported by BSF grants 9600210
and 2000352}

\begin{abstract}
We  study torus actions on symplectic manifolds with proper moment maps
in the case that each reduced space is two-dimensional.
We provide a complete set of invariants for such spaces.
\end{abstract}

\maketitle

\tableofcontents

\section{Introduction}
\labell{sec:intro}

Let a torus $T \cong (S^1)^{\dim T}$ act on a symplectic manifold
$(M,\omega)$ by symplectic transformations with moment map
$\Phi \colon M \to \ft^*$. We take the sign convention
\begin{equation} \labell{moment}
        \iota(\xi_M) \omega = - d \left< \Phi,\xi \right> .
\end{equation}
Here, $\xi$ is in the Lie algebra $\t$ of $T$
and $\xi_M$ is the vector field on $M$ induced by $\xi$.
Assume that the $T$-action is effective\footnote
{
 A group action is \emph{effective} if only the identity element 
 acts trivially.
}
on each connected component of $M$.
We call  $(M,\omega,\Phi)$
a \textbf{Hamiltonian $\mathbf T$-manifold}.\footnote
{
One sometimes allows $\omega$ to be degenerate.  Here we do not allow this.
}
If $\calT \subseteq \t^*$ is an open subset containing $\image \Phi$
and the map $\Phi \colon M \to \calT$ is proper,
then we call $(M,\omega,\Phi,\calT)$ a
\textbf{proper Hamiltonian $\mathbf{T}$-manifold}.

The \textbf{complexity} of $(M,\omega,\Phi)$ is
the difference  $k = \half \dim M - \dim T$; it  is half the dimension
of the reduced space  $\Phi\inv(\alpha)/T$
at a regular value $\alpha$ in $\image \Phi$.
For brevity, we call a complexity one proper Hamiltonian $T$-manifold
$(M,\omega,\Phi,\calT)$ a \textbf{complexity one space}
if $M$ is connected and $\calT$ is convex. 
A complexity one Hamiltonian $T$-manifold is \textbf{tall} 
if every reduced space is two dimensional, that is, if no reduced
space is a single point.

The simplest example of a complexity one space is a compact symplectic
surface $(\Sigma,\sigma)$ with no group action. The next
simplest example is the fiberwise circle action on a ruled surface.
More generally,
let $(M,\omega,\Phi)$ be a \emph{symplectic toric manifold}, that is,
a compact complexity zero Hamiltonian $T$-manifold. 
Let $P \to \Sigma$
be a principal $T$-bundle and let $\Theta \in \Omega^1(P,\t)$
be a connection one form.  For sufficiently large $k$, the form
$\tilde{\omega} = k \sigma + \omega + d \left< \Theta , \Phi \right>
 \in \Omega^2 (P \times M)$ descends to a symplectic form on $P \times_T M$
with moment map $\tilde{\Phi}([p,m]) = \Phi(m)$.  Then
$(P \times_T M , \tilde{\omega}, \tilde{\Phi} )$ is a tall complexity one
space.  For example, $\Sigma \times M$ is a tall complexity one space.
Finally, given a tall complexity one space, its equivariant
symplectic blow-up at any fixed point is also a tall complexity one space.

Symplectic toric manifolds are classified by their moment map images 
\cite{De}.  By Moser \cite{moser}, compact symplectic surfaces are 
classified by their genus and total area.
The next examples of complexity one spaces are
compact symplectic four manifolds with Hamiltonian  circle actions, which
were classified by the first author \cite{karshon:periodic},
following earlier work by Ahara, Hattori, and Audin
\cite{ah-ha,audin:paper,audin:book}.
In the algebraic category, complexity one actions (of possibly non-abelian
groups) were classified by Timash\"ev \cite{T1,T2}.
Chiang \cite{river} classified complexity one Hamiltonian actions 
of non-abelian groups on six-manifolds.
Li \cite{L} has obtained some classification results for certain Hamiltonian
circle actions on six manifolds.
However, a complete classification of
complexity two Hamiltonian torus actions would entail a classification
of four dimensional symplectic manifolds, which is not tractable.
See \cite{locun} for a more extensive list of related works.

Complexity one spaces are substantially more complicated than 
symplectic toric manifolds.
Symplectic toric manifolds provide a useful source of examples 
and counter-examples. We hope that complexity one
spaces will prove similarly useful, and that their greater complexity
$\begin{smallmatrix} \cdot\cdot \\ \smile \end{smallmatrix}$
will enable them to demonstrate phenomena that do not
occur on symplectic toric manifolds.
For example, all symplectic toric manifolds are K\"ahler.
However, there exist complexity one spaces with isolated fixed
points and no invariant K\"ahler structure \cite{T}.

This paper is the second in a series of papers in
which we study complexity one spaces.
Our goal is to classify these spaces. This consists of two parts.
First, uniqueness: 
we must determine whether or not two given spaces are equivariantly
symplectomorphic.
Second, existence: we must  provide a list of all complexity one spaces.

In \cite{locun} we obtained a \emph{local uniqueness} result:
we determined when two complexity one spaces
are equivariantly symplectomorphic over small subsets of $\t^*$.

In this paper we obtain a \emph{global uniqueness} result:
we provide invariants which determine when two tall complexity one
spaces are equivariantly symplectomorphic.

In our next paper of this series we will obtain \emph{global existence} 
results.  This will enable us to construct examples.

While the complexity one assumption is absolutely vital to our results,
the tall assumption is not.  In fact,  let $(M,\omega,\Phi,\calT)$
and $(M',\omega',\Phi',\calT)$ be any two complexity one spaces.
By removing every moment map fiber which consists of only one orbit,
we obtain tall complexity one spaces over an open subset of $\calT$.
We expect that the original manifolds will be equivariantly symplectomorphic
exactly if these tall complexity one spaces are equivariantly
symplectomorphic. 
However, 
the proof will require an additional ingredient: a variant
of Smale's theorem on the diffeomorphisms of $S^2$.

Whereas many complexity one spaces that one encounters
in nature are not tall, the ``tall" case is sufficient 
for constructing interesting examples.

\subsubsection*{Acknowledgement}
We would like to thank the referee for many useful suggestions
on the exposition.

\section{Statement and proof of ``Global Uniqueness"}
\labell{sec:statement and proof}

We  now describe the invariants of a tall complexity one space
$(M,\omega,\Phi,\calT)$.

An \textbf{isomorphism} between two Hamiltonian $T$-manifolds is 
an equivariant symplectomorphism that respects the moment maps.

Recall that Liouville measure on $M$ is given by integrating
the volume form $\omega^n / n!$ with respect to the symplectic orientation.
The \DH measure is the measure on $\calT$ obtained as the
push-forward of Liouville measure by the moment map.

The \textbf{stabilizer} of a point $x \in M$ is the closed subgroup
$\Stab(x) = \{\lambda \in T \mid \lambda \cdot x = x\}$.
The isotropy representation at $x$ is the linear representation
of $\Stab(x)$ on the tangent space $T_xM$.
Points in the same orbit have the same stabilizer,
and their isotropy representations are linearly symplectically isomorphic;
this isomorphism class is the \textbf{isotropy representation} of the
orbit.  An orbit is \textbf{exceptional} if every nearby orbit in the same
moment fiber $\Phi\inv(\alpha)$ has a strictly smaller stabilizer.  Let
$$M_\exc \subset M/T$$
denote the set of exceptional orbits.
The moment map induces a map $\Phibar \colon M_\exc \to \calT$
which is locally a proper embedding; this follows from the local
normal form theorem.

Let $M'_\exc$ denote the set of exceptional orbits of another tall
complexity one space.  An \textbf{isomorphism}
from $M_\exc$ to $M_\exc'$ is a homeomorphism that respects the moment maps
and sends each orbit to an orbit with the same isotropy representations.

\begin{Remark}
Assume, for simplicity, that $M$ is compact.
The orbit type decomposition of $M$ induces a stratification of $M_\exc$.
Given any stratum, $N$, the pre-image in $M$ of its closure $\ol{N}$
is a compact toric variety (with the action of $T/\Stab(N)$).
The moment map induces a diffeomorphism $\Phibar$ between $\ol{N}$
and the convex polytope $\Phibar(\ol{N}) \subseteq \t^*$.
If a stratum $N'$ is contained in $\ol{N}$ then $\Phibar(\ol{N'})$ is a face 
of $\Phibar(\ol{N})$.  So, topologically, $M_\exc$ is a union of
convex polytopes in $\t^*$, glued along faces, 
and $\Phibar \colon M_\exc \to \t^*$ restricts to the inclusion map
on each polytope.  This partially ordered collection of polytopes
is essentially equivalent to the notion of an X-ray as defined in~\cite{T}.
\end{Remark}

To define our other invariants, we need the following result.

\begin{Proposition} \labell{trivialize M mod T}
Let $(M,\omega,\Phi,\calT)$ be a tall complexity one space.
There exists a closed oriented surface $\Sigma$
and a map $f \colon M/T \to \Sigma$ so that
$$(\Phibar,f) \colon M/T \to (\image \Phi) \times \Sigma$$
is a homeomorphism, where $\Phibar$ is induced by the moment map.
Given two such maps $f$ and $f'$, there exists a homeomorphism
$\xi \colon \Sigma' \to \Sigma$ so that $f$ is homotopic to $\xi \circ f'$
through maps which induce homeomorphisms $M/T \to (\image \Phi) \times \Sigma$.
\end{Proposition}

The \textbf{genus} of the complexity one space is the genus of $\Sigma$.
Note that every  nonempty reduced space $\Phi\inv(\alpha)/T$
is homeomorphic to $\Sigma$.

A \textbf{painting} is a map $f \colon M_\exc \to \Sigma$ such that
the map
$$(\Phibar,f) \colon M_\exc \to \calT \times \Sigma$$
is one to one, where $M_\exc$ is the set of exceptional orbits.
Two paintings, $f \colon M_\exc \to \Sigma$ and 
$f' \colon M_\exc' \to \Sigma'$,
are \textbf{equivalent} if there exist an isomorphism
$i \colon M_\exc \to M'_\exc$ and a homeomorphism
$\xi \colon \Sigma \to \Sigma'$ such that the compositions
$\xi \circ f \colon M_\exc \to \Sigma'$ 
and $f' \circ i  \colon M_\exc \to \Sigma'$ are homotopic
\emph{through paintings}.

Proposition \ref{trivialize M mod T} implies that
there is a well-defined equivalence class of paintings associated to
every tall complexity one space;
just  restrict $f$ to $M_\exc$.

We can now state our main theorem:

\begin{Theorem} \labell{main-theorem}
Let $(M,\omega,\Phi,\calT)$ and $(M',\omega',\Phi',\calT)$
be tall complexity one spaces.
They are isomorphic if and only if they have the same
\DH measure, the same genus, and equivalent paintings.
\end{Theorem}

\begin{Remark}
A compact symplectic 4-manifold $M$ equipped with a Hamiltonian circle
action is tall if and only if the maximum and minimum of the moment map
are both attained on two dimensional surfaces, $\Sigma_\max$ and
$\Sigma_\min$.  Equivalently, it is tall exactly if it
can be obtained from a ruled surface by a sequence of
blowups  \cite{karshon:periodic}. 

By \cite{karshon:periodic}, 
the space is determined up to equivariant symplectomorphism
by the moment map values, genus, and symplectic areas of these surfaces,
together with the graph whose vertices correspond to isolated
fixed points in $M$ and are labeled by their moment map values 
and whose edges correspond to 2-spheres in $M$ with non-trivial
finite stabilizers and are labeled by the cardinalities of these stabilizers.

These invariants are equivalent to those of Theorem \ref{main-theorem}.
First, by Proposition \ref{trivialize M mod T},
the genus of $M$ is the genus of $\Sigma_\min$ and $\Sigma_\max$.
Second, $M$ and $M'$ have equivalent paintings exactly if $M_\exc$
is isomorphic to $M'_\exc$; this follows from the fact that $M_\exc$
is a union of intervals. By the Guillemin-Lerman-Sternberg formula
\cite[Section 3.5]{GLS}, 
the graph determines the \DH measure up to an affine function;
this function is determined by the moment map values
and symplectic areas of $\Sigma_\min$ and $\Sigma_\max$.
These, in turn, are determined by the \DH measure.
\end{Remark}

\begin{proof}[Proof of Proposition \ref{trivialize M mod T}]
By the convexity theorem for proper moment maps,
the image of $\Phi$ is convex \cite{LMTW}.
In \cite[Corollary 9.8]{locun}, we proved that
$\ol\Phi \colon M/T \to \image \Phi$ is topologically a
locally trivial bundle whose fiber is a closed oriented surface $\Sigma$.
Since the base is contractible and paracompact, the bundle is trivializable;
see \cite{husemoller}.

Any two trivializations, $(\ol\Phi,f)$ and $(\ol\Phi,f')$, differ by
a family of homeomorphisms $\xi_\beta \colon \Sigma' \to \Sigma$,
parametrized by $\beta \in \image \Phi$, that is determined by
$f(x) = \xi_{\Phibar(x)}(f'(x))$ for all $x$.
Pick any $\alpha \in \image \Phi$. Let
$f_t(x) = \xi_{(1-t)\Phibar(x)+t\alpha}(f'(x))$.
Then $f_0 = f$, and $f_1 = \xi_\alpha \circ f'$.
\end{proof}

We now prove our main theorem, using definitions and results
that we  develop later in the paper.

\begin{proof}[Proof of Theorem \ref{main-theorem}.]
By Proposition \ref{S to Q},
it is enough to show that the quotients $M/T$ and $M'/T$
are $\Phi$-diffeomorphic.
(See Definition \ref{Phi diffeo}.)

Now apply Proposition \ref{main prop}.
Let $N$ and $N'$ be painted surface bundles which are associated
to $M$ and $M'$. 
The quotients $M/T$ and $M'/T$ are $\Phi$-diffeomorphic
if and only if $N$ and $N'$ are isomorphic.
(See the definitions in Sections \ref{sec:S} and \ref{sec:exciting}.)

We can associate to $N$ and $N'$ \emph{smooth} paintings
$f \colon M_\exc \to \Sigma$ and $f' \colon M'_\exc \to \Sigma$,
and these paintings are equivalent.
Here we use Proposition \ref{bundle is painting} and
our assumption that $M$ and $M'$ have equivalent paintings.

This implies that $f$ and $f'$ are \emph{smoothly} equivalent,
by Proposition \ref{smooth to continuous}.
Then, $N$ and $N'$ are isomorphic, by Proposition \ref{bundle is painting}.
\end{proof}

\begin{Remark}
Fix a local homeomorphism $i \colon \calT \to \t^*$.  
Consider a symplectic manifold $(M,\omega)$
with a proper map $\Phi \colon M \to \calT$ such that $i \circ \Phi$
is a moment map for a $T$ action.  Suppose that $\Delta := \Phi(M)$
is contractible and the fibers $\Phi\inv(\alpha)$
are connected and two dimensional.  
Then $M/T$ is homeomorphic to $\Delta \times \Sigma$ 
as in Proposition \ref{trivialize M mod T},
and Theorem \ref{main-theorem} should remain true.
\end{Remark}


\section*{Part I:  Passing to $M/T$ }

\section{global structure of $M/T$.}

The goal of this paper is to determine
when two complexity one spaces are isomorphic.
In this section we  reduce this question to a simpler
question about their quotients, $M/T$ and $M'/T$.
To state this precisely, we introduce some
definitions from \cite{locun}.

\begin{Definition}
Let a torus $T$ act on oriented manifolds $M$ and $M'$ with
$T$-invariant maps $\Phi \colon M \to \t^*$ and $\Phi' \colon M' \to \t^*$.
A \textbf{$\bs \PhiT$-diffeomorphism} from $(M,\Phi)$ to  $(M,\Phi')$
is an orientation preserving equivariant diffeomorphism
$f \colon M \to M'$ that satisfies $\Phi' \circ f = \Phi$.
(Cf., \cite[Definition 3.1]{locun}.)
\end{Definition}

Let a compact torus $T$ act on a manifold $N$.
The quotient $N/T$ can be given
the quotient topology and a natural differential structure,
consisting of the sheaf of real-valued functions
whose pullbacks to $N$ are smooth.
We say that a map $h \colon N/T \to N'/T$ is \textbf{smooth}
if it pulls back smooth
functions to smooth functions; it is a \textbf{diffeomorphism}
if it is smooth and has a smooth inverse. See \cite{schwarz:IHES}.
If $N$ and $N'$ are oriented, the choice of an orientation on $T$
determines orientations on the smooth part of $N/T$ and $N'/T$.
Whether or not a diffeomorphism $f\colon N/T \to N'/T$ preserves orientation
is independent of this choice.

\begin{Definition} \labell{Phi diffeo}
Let $(M,\omega,\Phi,\calT)$ and $(M',\omega',\Phi',\calT)$
be complexity one Hamiltonian $T$-manifolds.
A $\mathbf{\Phi}$\textbf{-diffeomorphism} from $M/T$ to $M'/T$
is an orientation preserving diffeomorphism  $f \colon M/T \to M'/T$
such that ${\Phibar}' \circ f = \Phibar$ and such that
 $f$ and $f\inv$ lift to  $\PhiT$-diffeomorphisms
in a neighborhood of each exceptional orbit.
(Cf., \cite[Definition 4.1]{locun})
Here, $\Phibar$ and ${\Phibar}'$ are induced by the moment maps.
\end{Definition}

\begin{Proposition} \labell{S to Q}
Two tall complexity one spaces are isomorphic if and only if 
their quotients are $\Phi$-diffeomorphic
and their \DH measures are the same.
\end{Proposition}

The conditions are clearly necessary.
Proposition \ref{trivialize M mod T} implies that the restriction
$H^2(M/T,\Z) \to H^2(\Phi\inv(y)/T,\Z)$ is one to one for
every $y \in \image \Phi$.
Proposition \ref{S to Q} then follows from the results below,
which we proved in \cite[Propositions 3.3 and 4.2]{locun}:

\begin{Proposition}
Let $(M,\omega,\Phi,\calT)$ and $(M',\omega',\Phi',\calT)$
be complexity one spaces with the same \DH measure.
Assume  that the restriction map
$H^2(M/T,\Z) \to H^2(\Phi\inv(y)/T,\Z)$ is one to one for
some regular value $y$ of $\Phi$.
Then
\begin{enumerate}
\item
There exists an isomorphism from $M$ to $M'$
if and only if there exists a $\PhiT$-diffeomorphism from $M$ to $M'$.
\item
There exists an $\PhiT$-diffeomorphism from $M$ to $M'$
if and only if there exists a $\Phi$--diffeomorphism from $M/T$ to $M'/T$.
\end{enumerate}
\end{Proposition}

\section*{Part II: Abstract non-sense}

In \cite{locun}, we gave invariants that determine
the local pieces of a complexity one space.
In this paper we explain how these  pieces can be glued together.

This is analogous to classifying principle $G$-bundles
over a manifold $\calT$: locally they are trivial, and to determine
them globally, one needs to determine how the local pieces are glued together.
If $G$ is abelian, this is very easy:
the \v Cech cohomology, $\vH^i(\calT,G)$, is well defined for all
$i \geq 0$, and $G$-bundles are classified by $\vH^1(\calT,G)$.
Here, $G$ also denotes the sheaf of smooth functions to the group $G$.
When $G$ is not abelian, the $i$th \v{C}ech cohomology is only defined
for $i=0$ and $i=1$, and $G$-bundles are still classified
by $\vH^1(\calT,G)$.

A proper Hamiltonian $T$-manifold $(M,\omega,\Phi,\calT)$
determines a sheaf of non-abelian groups over $\calT$:
associate to an open subset $U \subseteq \calT$
the group of isomorphisms
of the preimage $\Phi\inv(U)$. The first cohomology of this sheaf
classifies Hamiltonian $T$-manifolds that are locally
isomorphic to $(M,\omega,\Phi,\calT)$, where ``locally"
means over small subsets of $\calT$.

We prefer, instead, to allow \emph{different} Hamiltonian $T$-manifolds
over each $U$, so that the isomorphisms between them form a \emph{groupoid}.
Besides being more elegant,
in that it does not single out one manifold above others,
this machinery lets us glue pieces of manifolds
without a-priori assuming that this can be done.
We use this in Section \ref{sec:globex}, where we prove
a technical result that will allow us, in subsequent papers,
to determine the full list of complexity one spaces
(``global existence").
In Sections \ref{sec:intro}--\ref{sec:smooth to continuous}, 
the reader may still choose to fix
a distinguished space and work with groups instead of groupoids.

We define sheaves of groupoids and their cohomology
in Sections \ref{sec:sheaves} and \ref{sec:cohomology}.
This straightforward extension of sheaves of abelian groups
and \v Cech cohomology sets up a convenient formalism.
These ideas are not new; closely related notions appear in
the literature.  See, for example, \cite[Chapter V]{Br}. 

In Sections \ref{sec:grommets}--\ref{sec:rigidify P to RP} 
we apply this formalism to a series of sheaves,
and show that they all have the same first cohomology.
This  reduces the classification of tall complexity one spaces
to a classification of simpler objects,
``painted surface bundles", which we classify directly
in Sections \ref{sec:S}--\ref{sec:smooth to continuous}.

\section{Sheaves of groupoids}
\labell{sec:sheaves}

A \textbf{groupoid} is a category $\cA$
where every arrow is invertible.
We let $\ob \cA$ denote the set of objects of $\cA$.
Given objects $A$ and $A'$,
let $\hom_\cA(A,A')$ denote the set of arrows with domain
$A$ and codomain $A'$,
and let $f \colon A \to A'$ denote an element of $\hom_\cA(A,A')$.
A \textbf{homomorphism} of groupoids is a functor.
Equivalently, given groupoids $\cA$ and $\cB$,
a homomorphism $\psi \colon \cA \to \cB$
consists of a map $\psi \colon \ob \cA \to \ob \cB$ and for
each  $A$ and $A'$ in $\ob \cA$
a map $\psi \colon \hom_\cA(A,A') \to \hom_\cB(\psi(A),\psi(A'))$
so that
$\psi(\id_A) = \id_{\psi(A)}$ and
$\psi(f) \circ \psi(f') = \psi(f \circ f')$.

\begin{Definition}
A \textbf{presheaf of groupoids} over a topological space $\calT$ assigns
\begin{enumerate}
\item  a groupoid $\cA(U)$ to every open subset $U \subseteq \calT$, and
\item  a homomorphism $\cA(\iota^V_U) \colon \cA(U) \to\cA(V)$,
called the \textbf{restriction map},
to every inclusion of open sets $V \subseteq U$.
\end{enumerate}
These must satisfy
\begin{enumerate}
\item
$\cA(\iota^U_U) = \id_{\cA(U)}$ for any open set $U$.
\item
$\cA(\iota^W_V) \circ \cA(\iota^V_U) = \cA(\iota^W_U)$
for any inclusions of open sets $W \subseteq V \subseteq U$.
\end{enumerate}
\end{Definition}

For an object $A$ and an arrow $f$ in $\cA(U)$,
let  $A|_V$ and $f|_V$ denote
$\cA(\iota^V_U)(A)$  and $\cA(\iota^V_U)(f)$, respectively.
Objects in $\cA(\calT)$ are called \textbf{global objects}.

\begin{Example} \labell{Gbundlesheaf}
Given a Lie group $G$,
define a presheaf  as follows:
the objects over $U \subset \calT$ are principle $G$ bundles over
$U$, and the arrows are  bundle isomorphisms.
Here, and in all other examples in this paper,
the restriction maps are given by restriction.
\end{Example}

\begin{Example} \labell{example:sheaves}
Let $T$ be a torus and $\calT \subseteq \t^*$ an open subset.
We may consider the following two presheaves.
In both, the objects over  $U \subseteq \calT$
are complexity one Hamiltonian $T$-manifolds
with proper moment maps to $U$.
The arrows in the first presheaf are equivariant symplectomorphisms
which respect the moment maps.
The arrows in the second presheaf are $\Phi$-diffeomorphisms
between the quotient spaces $M/T$.
\end{Example}

A sheaf is a presheaf where the arrows are determined by local data:

\begin{Definition} \labell{def:sheaf}
A \textbf{sheaf} over $\calT$ is a presheaf $\cA$
which satisfies the following two \textbf{sheaf axioms}.
Let $\{W_\alpha\}$ be a collection of open subsets of $\calT$.
Let $A$ and $A'$ be objects in $\cA(\cup W_\alpha)$.
\begin{enumerate}
\item
If $f \colon A \to A'$ and $g \colon A \to A'$ are arrows such that
$f|_{W_\alpha} = g|_{W_\alpha}$ for all $\alpha$, then $f = g$.
\item
Given a collection of arrows
$f_\alpha \colon A|_{W_\alpha} \to A'|_{W_\alpha}$
which are compatible on intersections,
there exists an arrow $f \colon A \to A'$
such that $f|_{W_\alpha} = f_\alpha$ for all $\alpha$.
\end{enumerate}
\end{Definition}

The second sheaf axiom says that arrows can be ``glued".
A sheaf $\cA$ has \textbf{gluable objects} if for every
collection $\{W_\alpha\}$ of open subsets of $\calT$,
objects $A_\alpha \in \cA(W_\alpha)$, and transition maps
$f_{\beta \alpha} \colon A_\alpha|_{W_\alpha \cap W_{\beta}}
\to A_\beta|_{W_\alpha \cap W_{\beta}}$, such that
$f_{\alpha \alpha} = \id_{A_\alpha}$ and
$f_{\gamma \beta} \circ f_{\beta \alpha} = f_{\gamma \alpha}$
on $W_\alpha \cap W_\beta \cap W_\gamma$, there exists an
object $A \in \cA(\cup W_\alpha)$ and isomorphisms
$\theta_\alpha \colon A|_{W_\alpha} \to A_\alpha$
for all $\alpha$ so that
$\theta_\beta \circ \theta_\alpha\inv= f_{\beta \alpha}$
on $W_\alpha \cap W_\beta$.

The  presheaves in Examples \ref{Gbundlesheaf} and
\ref{example:sheaves} are sheaves.
The first sheaf in  Example \ref{example:sheaves} has gluable
objects; the second sheaf does not.

\begin{Definition}
Let $\cA$ and $\cB$ be presheaves of groupoids over $\calT$.
A \textbf{map of presheaves} $\cF \colon \cA \to \cB$
assigns to every open subset $U \subseteq \calT$
a homomorphism $\cF(U) \colon \cA(U) \to \cB(U)$
so that for every inclusion of open sets $V \subseteq U$
we have $\cB(\iota^V_U) \circ \cF(U) = \cF(V) \circ \cA(\iota^V_U).$
\end{Definition}

\begin{Example}
Any equivariant symplectomorphism
  $M \to M'$ descends to a $\Phi$-diffeomorphisms
$M/T \to M'/T$. This gives a natural map from the first sheaf of
Example \ref{example:sheaves} to the second sheaf of that example.
\end{Example}

\section{\v Cech cohomology for sheaves of groupoids}
\labell{sec:cohomology}

In this section we define the first cohomology of a sheaf of groupoids
and give a geometric interpretation.

\begin{Definition}
Fix a sheaf of groupoids $\cA$ over $\calT$
and an open cover $\fU$ of $\calT$.

A \textbf{zero cochain} $a \in C^0(\fU,\cA)$ associates to each $U \in \fU$
an arrow $a_U$ in $\cA(U)$.
A \textbf{one cochain} $\alpha \in C^1(\fU,\cA)$ associates
to each $U \in \fU$ an object $A_U \in \cA(U)$
and to each pair $U,V \in \fU$ an arrow
$\alpha_{VU} \colon A_U|_{U\cap V} \to  A_V|_{U \cap V}$.

A \textbf{one cocycle} $\alpha \in Z^1(\fU,\cA)$ is a one cochain
that is \textbf{closed}, meaning that $\alpha_{UU} = \id_{A_U}$
for all $U \in \fU$ and that $\alpha_{WV} \circ \alpha_{VU} = \alpha_{WU}$
holds on $U \cap V \cap W$ for every triple $U, V, W \in \fU$,

The groupoid of zero cochains acts on the set of one cocycles by
$$ g_{VU} \mapsto
   \left( f_V|_{U \cap V} \right) \circ g_{VU}
   \circ \left( f_U|_{U \cap V} \right) \inv,$$
wherever this makes sense.
The \textbf{first cohomology} is the set of
equivalence classes under this action:
$$ \vH^1(\fU,\cA) := Z^1(\fU,\cA) / C^0(\fU,\cA). $$
\end{Definition}

An open cover $\fV$ is\footnote{
$\fV$ is \$$\backslash$\text{mathfrak V}\$.
}
a \textbf{refinement} of an open cover $\fU$ if every set $V \in \fV$
is a subset of a set $U \in \fU$.  As in the abelian case,
this induces a map in cohomology:

\begin{Lemma}
If $\fV$ is a refinement of $\fU$,
we get a well defined map in cohomology $\vH^1(\fU,\cA) \to \vH^1(\fV,\cA)$
for any presheaf of groupoids $\cA$.
\end{Lemma}

\begin{proof}
Choose any map $f \colon \fV \to \fU$ such that $V \subset f(V)$.
This map induces a map on one cocycles; given a cocycle $\alpha$,
simply restrict every object and
every map from $f(V)$ to $V$.
This clearly descends to a map on cohomology.
If $f' \colon \fV \to \fU$ is different map such that
$V \subset f'(V)$, the resulting cocycles  differ
by the following zero cochain:  associate to each $V \in \fV$
the restriction to $V$ of the arrow $\alpha_{f(V) f'(V)}.$
\end{proof}

Since the set of open covers is a directed set,
this makes $\vH^1(\fU,\cA)$ into a direct system of sets.
The \textbf{\v Cech cohomology} of $\calT$ with values in $\cA$,
denoted by $\vH^1(\calT,\cA)$, is defined to be the direct limit
of this direct system.

The following lemma is easy to check.

\begin{Lemma}
A map of sheaves of groupoids $f \colon \cA \to \cB$ induces a map of
cohomology $f_* \colon \vH^1(\calT,\cA) \to \vH^1(\calT,\cB)$.
\end{Lemma}

\begin{Lemma} \labell{classifies over fU}
Let $\cA$ be a sheaf of groupoids over $\calT$.
A global object in $\cA$ naturally determines
a class $[A] \in  \vH^1(\calT,\cA)$.
Two objects $A$ and $A'$ are isomorphic if and only if $[A] = [A']$.
If the sheaf has gluable objects, every class in  $\vH^1(\calT,\cA)$
arises in this way.
\end{Lemma}

\begin{proof}
A global object $A$ in $\cA(\calT)$ maps to the one cocycle
with object $A$ and arrow $\id_A$.
If
$f \colon A \to A'$ is an isomorphism, then $f$
is a zero cochain that interchanges the two cocycles.
On the other hand, if $A$ and $A'$
map to the same cohomology class. Then there exist isomorphisms
$f_U \colon A|_U \to A'|_U$, for $U$ in some open cover $\fU$,
so that $f_U|_{U \cap V} = f_V|_{U \cap V}$.
By the second sheaf axiom, this implies that $A$ is isomorphic to $A'$.
(See Definition \ref{def:sheaf}.)
Finally, the definition of gluable objects is chosen
exactly so that the map $A \mapsto [A]$ is onto for such sheaves.
\end{proof}

\section*{Part III: Sheaves of maps of $M/T$}


We wish to determine whether two tall complexity one spaces
are isomorphic.  By Proposition \ref{S to Q}, it is
enough to determine whether their quotients are $\Phi$-diffeomorphic.
By the results of \cite{locun}, the quotient $M/T$
is, topologically, a surface bundle over $\Phi(M)$.
If this were true in the $C^\infty$ category, it would be easy
to determine whether two quotients are $\Phi$-diffeomorphic.
Unfortunately, however, 
the quotient $M/T$ is naturally a manifold with corners
on the complement of the exceptional orbits, but not on the
exceptional orbits themselves. (See Lemma \ref{F smooth}.)

To overcome this difficulty, we
convert  our problem
to the problem of determining whether two
cohomology classes are  the same.
We define the sheaf $\cQ$ of $\Phi$-diffeomorphisms.
For  each open subset $U \subseteq \calT$, the objects in the groupoid $\cQ(U)$
are the tall complexity one proper Hamiltonian $T$-manifolds
over $U$; the arrows are $\Phi$-diffeomorphisms between their quotients.
By Lemma \ref{classifies over fU} and Proposition \ref{S to Q},
two tall complexity one space are isomorphic
if and only if they induce the same cohomology class
in $H^1(\calT,\cQ)$.

At first, this may not seem like a great improvement.
However, in this part of the paper,
we  gradually transform the sheaf $\cQ$
into one whose first cohomology we can compute.
To do this we  ``correct'' the smooth structure near the exceptional
orbits so that $M/T$ is a \emph{smooth} surface bundle.
(For more details, see Section \ref{sec:outline}.)
This process is unnatural; it relies on a choice of ``grommets".

\section{Grommets}
\labell{sec:grommets}

In this section, we define a new sheaf: 
$\Phi$-diffeomorphisms with grommets.

We fix an inner product on $\ft$, once and for all.
Let a closed subgroup $H \subseteq T$ act on $\C^n$
as a subgroup of $(S^1)^n$, with moment map
$\Phi_H  \colon \C^n \to \h^*$.
There exists an invariant symplectic form
on $$Y = T \times_H \C^n \times \h^0$$  with moment map
$$\Phi_Y ([t,z,\nu]) = \alpha + \Phi_H(z) + \nu,$$
where $\alpha \in \t^*$,
$\fh^0 \subseteq \ft^*$ is the annihilator of the Lie algebra $\fh$,
and we embed $\fh^*$ in $\ft^*$ using the metric.
The space $Y$ is called a \textbf{complexity one model}.

The \emph{local normal form theorem} \cite{GS:normal,marle} implies that
any orbit in a Hamiltonian $T$-manifold has a neighborhood
which is isomorphic to a neighborhood of the orbit $\{ [t,0,0] \}$
in some complexity one model $Y$.  This model is determined uniquely
up to permutations of the coordinates on $\C^n$; we call it
the \textbf{local model} associated to the orbit.

We recall the following definition from \cite[Definition 8.1]{locun}.

\begin{Definition}
\labell{def:grommet M}
Let $(M,\omega,\Phi)$ be a Hamiltonian $T$-manifold.
A \textbf{grommet} is a $\PhiT$-diffeomorphism $\psi \colon D \to M$
from an open subset $D$ of a local model $Y = T \times_H \C^n \times \h^0$
to an open subset of $M$.
\end{Definition}

Note that the domain $D$ need not contain the orbit $\{ [t,0,0] \}$.
Therefore, a grommet can be restricted to any open subset.

\begin{Definition}
\labell{def:grommeted M}
A complexity one Hamiltonian $T$-manifold
$(M,\omega,\Phi,\calT)$ is \textbf{grommetted}
if it is equipped with grommets whose images are disjoint and cover
the union of the exceptional orbits in $M$.
\end{Definition}

We think of this as attaching a grommet to the fabric of the
manifold at every exceptional orbit.  In real life, the fabric
can flow however it wants away from the grommets, but
at the grommet it can only spin;
this allows all the
necessary freedom of movement but prevents the fabric from ripping
at the points of stress.
Similarly, our grommets are designed to give enough freedom so
that we can approximate any map well, but still prevent us from having
to cope with the stress of really dealing with what happens at the
exceptional orbits.

We now define the  sheaf $\hcQ$ 
of $\Phi$-diffeomorphisms with grommets.
For each open set $U \subseteq \calT$
the objects in the groupoid $\hcQ(U)$
are the grommeted tall complexity one
proper Hamiltonian $T$-manifolds over $U$; the arrows are
$\Phi$-diffeomorphisms between their quotients
(which ignore their grommets).

The only difference between this sheaf $\hcQ$ and the sheaf $\cQ$ defined
in the beginning of Part III is that the objects in
$\hcQ$ carry grommets.
These sheaves have the same first cohomology.

\begin{Proposition}
\labell{lemma  hQ is Q}
The forgetful functor $\hcQ \to \cQ$
induces an isomorphism in cohomology,
$$\vH^1(\calT,\hcQ) \cong \vH^1(\calT,\cQ).$$
\end{Proposition}

The proof of this proposition 
uses an abstract sheaf-theoretic lemma:

\begin{Lemma} \labell{forgetful}
Let $i \colon \cA \to \cB$ be a map of sheaves such that:
\begin{enumerate}
\item For any open subset $U \subset \calT$ and objects
$A, A' \in \cA(U)$, $$i \colon \hom_\cA(A,A') \to \hom_\cB(i(A),i(A'))$$
is a bijection.
\item
For any open subset $U \subseteq \calT$ and
object $B \in \cB(U)$, every point in $U$
has  a neighborhood $V \subseteq U$ and an object
$A \in  \cA(V)$ so that $i(A)$ is isomorphic to $B|_V$.
\end{enumerate}

Then $i$ induces an isomorphism
$$ i_* \colon \vH^1(\calT,\cA) \to \vH^1(\calT,\cB).$$
\end{Lemma}

\begin{proof}
First we prove that $i_*$ is onto.
Let $\fU$ be a cover.
A cocycle $\beta \in  Z^1(\fU,\cB)$
associates to each $U \in \fU$ an object $B_U$ over $U$.
After passing to a refinement (which we still call $\fU$),
assumption (2) guarantees that  for each $U \in \fU$
there exists an object $A_U$ so that $i(A_U)$ is isomorphic to $B_U$.
By assumption (1), for each $U, V \in \fU$ there exists a unique
 $\alpha_{VU} \colon A_U|_{U \cap V} \to A_V|_{U \cap V}$
so that $i(\alpha_{VU}) = \beta_{VU}$.
Then $\alpha$ is a cocycle  and
$i(\alpha)$ is cohomologous to $\beta.$

Now we prove that $i_*$ is one-to-one.
It is enough to prove that the map 
$i_* \colon H^1(\fU,\cA) \to H^1(\fU,\cB)$ is one-to-one
for every cover $\fU$.
Suppose that $\alpha$ and $\alpha'$ are in $Z^1(\fU,\cA)$,
and that  $i(\alpha)$ and $i(\alpha')$ are cohomologous.
Then there exists a zero cochain which associates to
each $U \in \fU$ an arrow $b_U$ so that
$b_V \inv \circ i(\alpha_{VU}) \circ b_U = i(\alpha'_{VU})$
for all $V$ and $U$ in $\fU$.
By assumption (1), for each $U \in \fU$, there exists
a unique arrow $a_U$ such that $i(a_U) = b_U$.
Then $a_V \inv \circ \alpha_{VU} \circ a_U = \alpha'_{VU}$.
\end{proof}

\begin{proof}[Proof of Proposition \ref{lemma  hQ is Q}]
This follows from Lemma \ref{forgetful}
and from the fact that every complexity one proper Hamiltonian $T$-manifold
can be locally grommeted (see \cite[Lemma 8.4]{locun}).
\end{proof}

\section{$\Phi$-homeomorphisms}
\labell{sec:outline}

In this section, we list the sheaves that
we  need in this part of the paper. First, we must
consider maps which have the following form:

\begin{Definition} \labell{def:Phi homeo}
Let $M$ and $M'$ be complexity one Hamiltonian $T$-manifolds.
A $\mathbf{\Phi}$\textbf{-homeomorphism} between $M/T$ and $M'/T$
is a homeomorphism which sends each orbit to an 
orbit with the same local model, 
is a diffeomorphism off the set of exceptional orbits, 
respects the moment maps, and preserves the
orientation of the moment map fibers.
\end{Definition}

Let us now define the sheaf $\hcH$ of $\Phi$-homeomorphisms 
(with grommets).
For each open subset $U \subseteq \calT$, the objects
in the groupoid $\hcH(U)$ are the grommeted tall complexity one proper
Hamiltonian $T$-manifolds over $U$; the arrows are $\Phi$-homeomorphisms
between their quotients (which ignore the grommets).

We  work with a sequence of subsheaves of $\hcH$.
Each subsheaf has  same objects as  $\hcH$ does, but
the arrows  are $\Phi$-homeomorphisms which satisfy additional conditions.
The first cohomology of each sheaf is isomorphic to
that of $\cQ$.

We list the names and symbols for the sheaf $\hcH$
and its relevant subsheaves: \medskip

\begin{tabular}{lc}
$\Phi$-homeomorphisms (with grommets) & $\hcH$ \medskip \\ 
$\Phi$-diffeomorphisms (with grommets) & $\hcQ$   \\
locally rigid $\Phi$-homeomorphisms \quad & $\cRQ$   \\
local stretch maps & $\cE$ \\
locally sb-rigid  $\Phi$-homeomorphisms & $\cRP$   \\
sb-diffeomorphisms & $\hcP$
\end{tabular}\medskip

In Section \ref{sec:grommets} we defined the sheaf $\hcQ$ of 
$\Phi$-diffeomorphisms with grommets; note that it is 
a subsheaf of the sheaf of $\Phi$-homeomorphisms (with grommets).
We next  restrict to the subsheaf $\cRQ \subset \hcQ$
consisting of those $\Phi$-homeomorphisms that are, roughly
speaking,  given by
``rigid rotations" along the exceptional orbits.
Then we  extend to the sheaf $\cE$ of maps that are given by
``rotations and stretches" along the exceptional orbits.
We then restrict again to a sheaf $\cRP$ of ``rigid rotations",
except that these are defined differently, in such a way that they
are smooth with respect to a new differential structure
that makes $M/T$ into a smooth manifold with corners.
Finally, we extend to the sheaf $\hcP$ of \emph{all} maps that are smooth
with respect to the new differential structure.

We define the rest of these sheaves in Sections
\ref{sec:def rigid}--\ref{sec:pi diffeo}.
We also show that we have  the inclusions
\begin{equation} \labell{inclusions of sheaves}
\hcQ \supseteq \cRQ \subseteq \cE \supseteq \cRP \subseteq \hcP.
\end{equation}
In Sections \ref{sec:rigidify Q to RQ}--\ref{sec:rigidify P to RP}
we show that each of these inclusions induces an isomorphism on
the first cohomology. 

\section{Locally rigid $\Phi$-homeomorphisms}
\labell{sec:def rigid}

In this section we define the second sheaf in our sequence
of sheaves: locally rigid $\Phi$-homeomorphisms.

Let an $h$ dimensional closed subgroup $H \subseteq T$ act on $\C^{h+1}$
through an inclusion map
$$\rho = (\rho_0,\ldots,\rho_h) \colon H \to (S^1)^{h+1}$$
so that the origin is an exceptional orbit.
Let $R_\rho \subset U(h+1)$ denote the group of
unitary transformations that commute with the $H$ action.
Up to permutation there are only two
possibilities: either the $\rho_i$ are all different, or they
are all different except that $\rho_0 = \rho_1$.
In the former case, 
$$ R_\rho = (S^1)^{h+1} \subset U(h+1). $$
In the latter case, 
$$ R_\rho = U(2) \times (S^1)^{h-1} \subset U(h+1). $$
The group
$$ R_Y := T \times_H R_\rho $$
acts on the complexity one model
$ Y := T \times_H \C^{h+1} \times \h^0.$
Define $$\Rbar = R_Y/T = R_\rho/H.$$  The action of $R_Y$ on $Y$
descends to an action of $\Rbar$ on $Y/T$ through the
short exact sequence
$ 1 \to T \to R_Y \to \Rbar \to 1.$

\begin{Remark}
In \emph{tall} complexity one models, the $\rho_i$ are always
different, so that $R_\rho = (S^1)^{h+1}$ and $\Rbar \cong S^1$.
However, in this section we allow the general case 
because it does not require much extra work and will be useful in
subsequent papers.
\end{Remark}

We would like to use $\Rbar$ to define 
``locally rigid $\Phi$-homeomorphisms"
between grommeted complexity one
Hamiltonian $T$-manifolds.  However, even if the manifolds are 
isomorphic, the domains of the
given grommets might lie in different models, whereas the elements
of $R_Y$ are maps from a model to itself.
We solve this problem by
passing to \emph{sub-grommets}.  We define these in the next several
pages.

\begin{Lemma} \labell{standard}
Let $Y = T \times_H \C^{h+1} \times \fh^0$ be a complexity one model.
Let $E = \{ [t,y,\mu] \}$ be an exceptional orbit,
where $y_i = 0$  exactly if $0 \leq i \leq k$.

The local model associated to $E$ is
$Y_E = T \times_K \C^{k+1} \times \fk^0$,
where
\begin{equation} \labell{K}
 K = H \cap ((S^1)^{k+1} \times \{1\}^{h-k})
\end{equation}
acts on $\C^{k+1}$ as  the restriction
of the $H$ action on $\C^{h+1}$ to the first $k+1$ coordinates.
\end{Lemma}

\begin{proof}
Clearly, $K$ is the stabilizer of $E$.
Let $\cO = H \cdot y$.  Then $T_y \cO = \{0\}^{k+1} \times V$,
where $V$ is an $(h- \dim K)$-dimensional isotropic subspace of $\C^{h-k}$.
The symplectic slice is
$$(T_y\cO)^\omega /T_y \cO \cong \C^{k + 1}  \times \C^m,$$
where $m = \dim K - k$ and $K$ acts trivially on $\C^m$.
If the orbit $E$ is exceptional, then $m = 0$.
\end{proof}

\begin{Definition} \labell{def:canonical inclusion}
Let
$Y = T \times_H \C^{h+1} \times \fh^0$
be a complexity one model.  Let
$Y_E = T \times_K \C^{k+1} \times \fk^0$
be the local model associated to an exceptional orbit
$E = \{ [t,y,\mu] \} \subset Y,$
where $y_i = 0$ exactly if $0 \leq i \leq k$.

A \textbf{canonical inclusion} is a $\PhiT$-diffeomorphism
$$\Lambda \colon D_E \to
   Y = T \times_H \left( \C^{k+1} \times \C^{h-k} \right) \times \fh^0$$
such that
\begin{equation}\labell{Lambda}
\Lambda \left( [t,z,\nu] \right)
=  \left[ t, z, f([t,z,\nu]), \alpha([t,z,\nu]) \right],
\end{equation}
where  $D_E \subseteq Y_E$ is an open subset, and
$f \colon D_E \to  \R_{>0}^{h-k}$ and $\alpha \colon D_E \to \fh^0$
are $R_{Y_E}$ invariant functions.
\end{Definition}

The composition of two canonical inclusions is a canonical inclusion.

\begin{Lemma} \labell{canonical inclusion}
Let $Y$ be  a complexity one model, and let $Y_E$ be the local model
associated to an exceptional orbit $E \subset Y$.

There exists a
a canonical inclusion $\Lambda \colon D_E \to Y$
on some neighborhood $D_E$ of $\{[t,0,0]\}$ in $Y_E$.
On any open subset $D_E$ of $Y_E$ there exists at most one canonical
inclusion $\Lambda \colon D_E \to Y$.
\end{Lemma}

\begin{proof}
Here, we use the notation of  Definition \ref{def:canonical inclusion}.

Let $\eta_j \in \fh^*$ denote the weights for the $H$ action on $\C^{h+1}$.
A moment map for $Y$ is
\begin{equation} \labell{mm on Y}
\Phi_Y([t,w,\alpha]) = \half \sum_{j=0}^h \eta_j |w_j|^2 + \alpha.
\end{equation}

The weights for the $K$ action on $\C^{k+1}$
are $\iota^* \eta_j$, for $0 \leq j \leq k$, where $\iota^*$
is the dual to the inclusion map $\iota \colon \fk \to \fh$.
The corresponding moment map for the local model $Y_E$ is
\begin{equation} \labell{mm on YE}
\begin{split}
\displaystyle{ \Phi_E([t,z,\nu])
   \ = \
\Phi_Y(E) }
  &+  \displaystyle{ \half \sum_{j=0}^k \iota^* \eta_j |z_j|^2 + \nu. } \\
   =
\displaystyle{ \half \sum_{j=k+1}^h \eta_j |y_j|^2 + \mu }
  &+  \displaystyle{ \half \sum_{j=0}^k \iota^* \eta_j |z_j|^2 + \nu. }
\end{split}
\end{equation}

Given \emph{any} pair of smooth $T$ invariant functions,
$f \colon D_E \to \R_{> 0}^{h-k}$
and $\alpha \colon D_E \to \h^0$, we can define
a smooth $T$-equivariant map $\Lambda \colon D_E \to Y$  by
the formula \eqref{Lambda}.
This map satisfies $ \Phi_Y \circ \Lambda  = \Phi_E$ if and only if
\begin{equation} \labell{same mm first}
\begin{split}
 & \half \left( \sum_{j=0}^k \eta_j |z_j|^2 +
             \sum_{j=k+1}^h \eta_j f_j^2 \right) + \alpha \\
 & \quad \ = \half \sum_{j=k+1}^h \eta_j |y_j|^2 + \mu
 + \half \sum_{j=1}^k \iota^* \eta_j |z_j|^2 + \nu.
\end{split}
\end{equation}

The metric induces a decomposition
$$ \ft^* = \fk^* \oplus \left( \fh^* \cap \fk^0 \right) \oplus \fh^0.$$

By Lemma \ref{standard},
the weights $\eta_j$, for $k < j \leq h$, lie in $\fh^* \cap \fk^0$.
Hence, the $\k^*$ components of the left and right hand sides
of \eqref{same mm first} automatically agree; they are both equal to
$\half \sum_{j=0}^k \iota^* \eta_j |z_j|^2$
Therefore, Equation \eqref{same mm first} is equivalent to the equations
\begin{equation} \labell{same mm}
\half \sum_{j=k+1}^h \eta_j f_j^2 =
\half \sum_{j=k+1}^h \eta_j |y_j|^2
+ \half \sum_{j=0}^k (\iota^* \eta_j - \eta_j) |z_j|^2 + \pi(\mu + \nu)
\end{equation}
in $\fh^* \cap \fk^0$ and
$$ \alpha = \mu + \nu - \pi(\mu + \nu)$$
in $\fh^0$,
where $\pi$ is the projection from $\fk^0$ to $\fk^0 \cap \fh^*$.

We find $f_i^2$ by solving the system \eqref{same mm} of linear equations.
The solution exists and is unique because the coefficient vectors
$\eta_j$, for $j=k+1,\ldots,h$, are a basis of $\fh^* \cap \fk^0$.
Since $f_j^2 = |y_j|^2 > 0$ when $z=0$ and $\nu=0$,
we can take smooth positive  square roots of $f_j^2$ near $[t,0,0]$.
Finally, the functions $f_j$ are $R_{Y_E}$-invariant
because the equation \eqref{same mm} is invariant and the solution
is unique.

It is clear that the resulting map $\Lambda$
preserves the orientation on each fiber.
To show that $\Lambda$ is a diffeomorphism,
it is enough to show that $\Lambda$ is a submersion.

Consider the map $H \to (S^1)^{h-k}$ obtained by the inclusion into
$(S^1)^{h+1}$ followed by the projection to the last $h-k$ coordinates.
The kernel of this map is $K$. Since $\dim K = k$, the map must be onto.
That is, we have a short exact sequence
$$ 1 \to K \to H \to (S^1)^{h-k} \to 1 .$$
This implies that the natural inclusion $\R_+^{h-k} \to (\Cc)^{h-k}$
gives rise to a $\PhiT$-diffeomorphism
$$ Y \supseteq 
T \times_H \left( \C^{k+1} \times (\Cc)^{h-k} \right) \times \fh^0
\cong T \times_K \C^{k+1} \times \R_+^{h-k} \times \fh^0.$$

In these coordinates, $\Lambda$ has the form
$[t,z,\nu] \mapsto [t,z,f,\mu + \nu - \pi(\mu + \nu)]$.
It is easy to check that this is a submersion.
\end{proof}

\begin{Definition} \labell{def:compatible grommet}
Let $M$ be a grommeted complexity one Hamiltonian $T$-manifold.
Let $\Lambda \colon D_E \to Y$ be a canonical inclusion whose image
is contained in the domain of a grommet $\psi \colon D \to M$.
Define $\psi_E  \colon D_E \to M$
by $\psi_E = \psi \circ \Lambda$.
We call the induced map $\psibar_E \colon D_E/T \to M/T$
a \textbf{sub-grommet}.
\end{Definition}

\begin{Definition} \labell{def:rigid on quotient}
Let $M$ and $M'$ be grommeted complexity one
Hamiltonian $T$-manifolds.  A $\Phi$-homeomorphism
$ f \colon M/T \to M'/T $
is \textbf{locally rigid}
if for every exceptional orbit $E \in M/T$ and any pair of sub-grommets
$$ \ol\psi \colon D/T \to M/T \quad \text{ and } \quad
   \ol\psi' \colon D/T \to M'/T $$
whose images contain $E$ and $f(E)$,
there exists a  smooth function $R \colon \t^* \to \Rbar$
 such that
 \begin{equation} \labell{rigid}
 (\ol\psi')\inv \circ f \circ \ol\psi(y)
 =  R(\Phibar_Y(y)) \cdot y
 \end{equation}
on some neighborhood of $\ol\psi\inv(E)$.
Here, both sub-grommets have the same domain $D/T \subset Y/T$,
and $\Phibar_Y$ is induced by the moment map on $Y$.
\end{Definition}

\begin{Remark}
By Lemma \ref{canonical inclusion}, for any exceptional orbit $E$
we can always find sub-grommets
$\ol\psi \colon D/T \to M/T$ and $\ol\psi' \colon D/T \to M'/T$,
with the same domain,
whose images contain $E$ and $f(E)$.
Since a canonical inclusion induces an $\Rbar$-equivariant map
on quotients,
Equation \eqref{rigid}  holds
for either every such $\ol\psi$ and $\ol\psi'$
or for no such $\ol\psi$ and $\ol\psi'$.
\end{Remark}

The following lemma is straightforward.

\begin{Lemma}
Every locally rigid $\Phi$-homeomorphism is a $\Phi$-diffeomorphism.
\end{Lemma}

We now define the sheaf $\cRQ$ of locally rigid $\Phi$-homeomorphisms.
For each open subset $U \subseteq \calT$, the objects in the
the groupoid $\cRQ(U)$ are the grommeted
tall complexity one proper Hamiltonian $T$-manifolds over $U$;
the arrows are the locally rigid $\Phi$-homeomorphisms between 
their quotients.

Our main claim, which we prove in Section \ref{sec:rigidify Q to RQ}, is

\begin{Proposition} \labell{lemma RQ is hQ}
The inclusion $\cRQ \subset \hcQ$ induces an isomorphism
$$ \vH^1(\calT,\cRQ) \cong \vH^1(\calT,\hcQ).$$
\end{Proposition}

\section{Local stretch maps}
\labell{sec:stretch}

In this section we define the third sheaf in our sequence:
local stretch maps.

Consider a tall complexity one model
$ Y = T \times_H \C^{h+1} \times \h^0 $
with moment map
$ \Phi_Y.$

\begin{Lemma}\labell{lem:P exact}
There exists a unique monomial $P \colon \C^{h+1} \to \C$,
called the \textbf{defining monomial}, such that
\begin{equation} \labell{P}
P(z) = \prod_{j=0}^h z_j^{\xi_j},
\end{equation}
with $\xi_j \geq 0$ for all $j$, such that the following sequence is
exact:
\begin{equation} \labell{P exact}
 1 \to H \stackrel{\rho}{\hookrightarrow} (S^1)^{h+1}
         \stackrel{P}{\to} S^1 \to 1.
\end{equation}
\end{Lemma}

\begin{proof}
See \cite[Lemma 5.8]{locun}.
\end{proof}

The \textbf{trivializing homeomorphism} is the map
$$ F = (\Phibar_Y, \ol{P}) \colon Y/T \to (\image \Phi_Y) \times \C ,$$
where $\Phibar_Y$ is induced from the moment map and $\ol{P}$
is induced from the defining monomial.
(See \cite[Definitions 5.12 and 6.4]{locun}.)

\begin{Lemma} \labell{F smooth}
The trivializing homeomorphism $F$ is a homeomorphism,
and it is a diffeomorphism off the set of exceptional orbits.
\end{Lemma}

\begin{proof}
See \cite[Lemmas 6.2 and 7.1]{locun}.
\end{proof}

A grommet on a tall complexity one Hamiltonian $T$-manifold induces
a ``coordinate chart" on $M/T$.

\begin{Definition} \labell{def:pi compatible}
Let $\psi \colon D \to M$ be a grommet
with $D \subseteq Y$, and let $F \colon Y/T \to \image \Phi_Y \times \C$
be the trivializing homeomorphism.
Let $B = F(D) \subset \image \Phi_Y \times \C$ 
and define $\varphi \colon B \to M/T$ by
$$ \varphi := \psibar \circ F\inv.$$ 
Then $\varphi$ is a homeomorphism  onto an open subset of $M/T$;
it is the \textbf{surface bundle grommet}
associated to $\psi$.
\end{Definition}

We can now give our main definition.

\begin{Definition} \labell{def:stretch map}
Let $M$ and $M'$ be grommeted tall complexity one
Hamiltonian $T$-manifolds.
A $\Phi$-homeomorphism $f \colon M/T \to M'/T$
is a \textbf{local stretch map} if for every exceptional
orbit $E \in M/T$ and the pair of associated surface bundle grommets
$\varphi \colon B \to M/T$ and $\varphi' \colon B' \to M'/T$
whose images contain $E$ and $f(E)$,
there exists a function $R \colon \ft^* \to S^1$
and an $S^1$-invariant function $\lambda \colon \t^* \times \C \to \R_{>0}$
such that
\begin{equation} \labell{form of stretch map}
 (\varphi')\inv \circ f \circ \varphi (\alpha,z)
 = (\alpha , R(\alpha) \lambda(\alpha,z) \cdot z),
\end{equation}
on some neighborhood of $\varphi\inv(E) \subseteq \ft^* \times \C$.
\end{Definition}

We need the following lemma.

\begin{Lemma} \labell{rigid is stretch a}
Let $\Lambda \colon D_E \to Y$ be a canonical inclusion 
for $D_E \subset Y_E$.
Let $F \colon Y/T \to \image \Phi_Y \times \C$ 
and $F_E  \colon Y_E/T \to \image \Phi_{Y_E} \times \C$ 
be the trivializing homeomorphisms.
Then there exists an $S^1$ invariant function 
$\lambda \colon D_E \to \R_{>0}$ such that
$$ F \circ \ol{\Lambda} \circ F_E\inv(\alpha,z) 
   = (\alpha, \lambda(\alpha,z) \cdot z),$$
for all $(\alpha,z) \in  F_E(D_E/T) \subset \t^* \times \C$.
\end{Lemma}

\begin{proof}
By Lemma \ref{standard},  the defining
monomial $P_E$ for $Y_E$ consists of the first $k+1$ factors
the the defining monomial $P$ for $Y$.
Hence, for $u=[t,z,\nu] \in Y_E$, we have
$$ P_E(u) = P_E([t,z,\nu]) = \prod_{j=0}^k z_j^{\xi_j} $$
and
$$ P( \Lambda(u)) = P([t,z,f,\tilde{\mu}])
   = \prod_{j=0}^k z_j^{\xi_j} \prod_{j=k+1}^h f_j^{\xi_j} $$
for some $\xi_0,\ldots,\xi_h$.
Thus, the lemma holds  with 
$$ \lambda(\alpha,z) 
   = \prod_{j=k+1}^h f_j({F_E}\inv(\alpha,z))^{\xi_j} .$$
\end{proof}

Carefully unwinding the definitions, this leads 
to the following lemma.

\begin{Lemma} \labell{rigid is stretch b}
Every locally rigid $\Phi$-homeomorphism is a local stretch map.
\end{Lemma}

We introduce the sheaf $\cE$ of local stretch maps.
For each open subset $U \subseteq \calT$,
the objects in the groupoid
$\cE(U)$ are the grommeted tall
complexity one proper Hamiltonian $T$-manifolds $M_U$;
the arrows are the local stretch maps between their quotients.

Our main claim, which we prove in Section \ref{sec:rigidify E}, is

\begin{Proposition} \labell{lemma RQ is E}
The inclusion $\cRQ \subset \cE$ induces an isomorphism
$$ \vH^1(\calT,\cRQ) \cong \vH^1(\calT,\cE).$$
\end{Proposition}

\section{Locally sb-rigid $\Phi$-homeomorphisms}
\labell{sec:sb-rigid}

In this section we define the fourth sheaf in our sequence: 
locally sb-rigid $\Phi$-homeomorphisms.
Here, ``sb" stands for ``surface bundle".

\begin{Definition} \labell{def:sb rigid}
Let $M$ and $M'$ be grommeted tall complexity one proper Hamiltonian
$T$-manifolds.
A $\Phi$-homeomorphism  $f \colon M/T \to M'/T$
is  \textbf{locally sb-rigid} if for every exceptional orbit $E \in M/T$ 
and the pair of associated surface bundle grommets
$\varphi \colon B \to M/T $ and $\varphi' \colon B' \to M'/T$
whose images contain $E$ and $f(E)$ (see Definition \ref{def:pi compatible}),
there exists a smooth function $R \colon \ft^* \to S^1$ such that
\begin{equation} \labell{composition}
 (\varphi')\inv \circ f \circ \varphi (\alpha,z)
 = (\alpha , R(\alpha) \cdot z)
\end{equation}
on some neighborhood of $\varphi\inv(E) \subseteq \ft^* \times \C$.
\end{Definition}

From this and Definition \ref{def:stretch map}
we immediately get the following result:

\begin{Lemma}
Every locally sb-rigid $\Phi$-homeomorphism is a local stretch map.
\end{Lemma}

We now define the sheaf $\cRP$ of locally sb-rigid $\Phi$-homeomorphisms.
For each open set $U \subseteq \calT$,
the objects in the groupoid $\cRP(U)$ are the grommeted tall complexity one
proper Hamiltonian $T$-manifolds over $U$; the arrows are
the locally sb-rigid $\Phi$-homeomorphisms  between their quotients.

Our main claim, which we prove in Section \ref{sec:rigidify E}, is

\begin{Proposition} \labell{lemma RP is E}
The inclusion $\cRP \to \cE$ induces an isomorphism
$$ \vH^1(\calT,\cRP) \cong \vH^1(\calT,\cE) .$$
\end{Proposition}

\section{sb-diffeomorphisms}
\labell{sec:pi diffeo}

In this section we define the fifth sheaf in our sequence: 
sb-diffeomorphisms. Again, ``sb" stands for ``surface bundle".

\begin{Definition} \labell{def:pi diffeo}
Let $M$ and $M'$ be grommeted tall complexity one  Hamiltonian
$T$-manifolds.
A $\Phi$-homeomorphism  $f \colon M/T \to M'/T$
is an  \textbf{sb-diffeomorphism} if for every
pair of associated surface bundle grommets
$\varphi \colon B \to M/T$ and $\varphi' \colon B' \to M'/T$,
the composition
$$ g = (\varphi')\inv\circ f \circ \varphi $$
is a diffeomorphism.
\end{Definition}

The definition clearly implies
\begin{Lemma}
If a $\Phi$-homeomorphism is locally sb-rigid 
then it is an sb-diffeomorphism.
\end{Lemma}

We now define the sheaf $\hcP$ of sb-diffeomorphisms.
For each open subset $U \subseteq \calT$,
the objects in the groupoid $\hcP(U)$  are the
grommeted tall complexity one proper Hamiltonian $T$-manifolds over $U$;
the arrows are the sb-diffeomorphisms between their quotients.

Our main claim, which we prove in Section \ref{sec:rigidify P to RP}, is

\begin{Proposition} \labell{lemma RP is hP}
The inclusion $\cRP \subset \hcP$ induces an isomorphism
$$ \vH^1(\calT,\cRP) \cong \vH^1(\calT,\hcP).$$
\end{Proposition}

\section*{Part IV: Rigidification}

In this evil part of the paper we prove the main propositions
stated in the previous part: Propositions
\ref{lemma RQ is hQ}, \ref{lemma RQ is E}, \ref{lemma RP is E},
and \ref{lemma RP is hP}.
In Section \ref{sec:reducing} we prove a technical lemma,
which we use four times, in Sections \ref{sec:rigidify Q to RQ},
\ref{sec:rigidify E}, and \ref{sec:rigidify P to RP},
to show that each of the inclusions of sheaves
$$ \hcQ \supseteq \cRQ \subseteq \cE \supseteq \cRP \subseteq \hcP $$
induces an isomorphism on $H^1$.  
The reader may choose to skip to the next part of the paper
(Section \ref{sec:S}) on first reading.

Checking the assumptions of the technical lemma
involves \emph{isotopies} of maps on complexity one quotients.
In previous sections we defined what it means for a map between
complexity one quotients to be a $\Phi$-homeomorphism,
a $\Phi$-diffeomorphism, a locally rigid $\Phi$-homeomorphism,
a local stretch map, an sb-diffeomorphism, or a locally sb-rigid 
$\Phi$-homeomorphism.
In Sections \ref{sec:rigidify Q to RQ}--\ref{sec:rigidify P to RP}
we will need to extend these notions to ``isotopies" and to being
``rigid at a point".  Also, we will need to have similar notions
on ``surface bundle models". We start with the most general definition:

\begin{Definition} \labell{def:isotopy Phi homeo}
Let $(M,\omega,\Phi,\calT)$ and $(M',\omega',\Phi',\calT)$
be complexity one Hamiltonian $T$-manifolds.  
A map $F$ from an open subset of $[0,1] \times M/T$
to an open subset of $M'/T$ is an
\textbf{isotopy of} {\boldmath$\Phi$}\textbf{-homeomorphisms}
if $f_t(\cdot) = F(t,\cdot)$ is a $\Phi$-homeomorphism for each $t$
and $F$ is smooth on the complement of the exceptional orbits
(as a map between manifolds with corners).
\end{Definition}

\begin{Remark} \labell{single Phi homeo}
We may consider a single $\Phi$-homeomorphism as an isotopy
which is independent of the parameter $t$.
\end{Remark}

Note that if $f_t(\cdot) = F(t,\cdot)$ is an isotopy
of $\Phi$-homeomorphisms from $M/T$ to $M'/T$ 
and $E$ is an exceptional orbit in $M$
then the exceptional orbit $f_t(E)$ 
remains the same as $t$ varies continuously.

\section{Reducing from sheaf to sheaf}
\labell{sec:reducing}

In this section we prove a lemma
which guarantees  that  the first cohomology of two different
sheaves agree if the sheaves satisfy
certain technical conditions.
Let $\cA$ and $\cB$ be sheaves of groupoids over $\calT$ with the same
objects and with the arrows in $\cA$ forming subsets of the arrows
in $\cB$.
The inclusion maps $\cA(U) \hookrightarrow \cB(U)$
induce a map from $\vH^1(\fU,\cA)$ to $\vH^1(\fU,\cB)$.
We describe a condition which guarantees that this map is a bijection.

Let us begin with a simple  analogue for sheaves of abelian groups.
Consider a manifold $M$ and two abelian Lie groups $H \subset G$.
Denote the sheaves of smooth functions to $H$ and to $G$ by
$\bfH$ and $\bfG$, respectively.
Suppose that  $H$ is a smooth deformation retract of $G$;
for example, $H=S^1$ and $G=\C^*$.  Then the natural map between the
cohomology groups $\vH^i(M,\bfH)$ and $\vH^i(M,\bfG)$
is an isomorphism for all $i \geq 1$.

One proof relies on the the fact that these sheaves satisfy the following
condition: if we are given an open set $Z$,
a smooth function $\beta \colon Z \to G$, and a
pair of open sets $X$ and $Y$ such that $\overline{X} \cap \overline{Y} =
\emptyset$,
we can define a new smooth function $\beta' \colon Z  \to G$ which satisfies
the following conditions:
\begin{enumerate}
\item
For all $x \in X \cap Z$,  $\beta'(x)$  is in $H$.
\item
For all $y \in Y \cap Z$,  $\beta'(y) = \beta(y)$.
\item
For all $u \in Z$,
if  $\beta(u)$ is in $H$,  then  $\beta'(u)$ is also in $H$.
\end{enumerate}

This fact is easy to show:
By assumption, there exists a smooth map 
$F \colon G \times [0,1] \to G$ so that
$F(g,0) = g$ for all $g \in G$, that $F(g,1) \in H$ for all $g \in G$,
and that $F(h,t) = h$ for all $h \in H$ and $t \in [0,1]$.
Since $\overline{X} \cap \overline{Y} = \emptyset$,
we can find a smooth function $\lambda \colon M \to [0,1]$
such that $\lambda(x) = 1$ for all $x \in X$, and $\lambda(y) = 0$
for all $y \in Y$.
Now, simply define $\beta'(u) = F(\beta(u),\lambda(u))$.

Our main technical lemma, which appears below, is a generalization of
the fact that this  condition is itself enough to prove that for
all $i \geq 1$ the cohomologies $\vH^i(M,\bfG)$ and $\vH^i(M,\bfH)$ 
are isomorphic.

\begin{Lemma} \labell{criterion for tech lemma}
Let $\cA$ and $\cB$ be sheaves of groupoids on $\calT$
with the same objects and such that  the arrows in $\cA$ are subsets
of the arrows in $\cB$.

Suppose that for every open cover $\fU = \{ U_i \}$ and cocycle
$\beta \in \vZ^1( \fU , \cB)$ there exists an open cover $\{U_i'\}$
such that $\ol{U_i'} \subset U_i$ for all $i$ and such that
the following holds:

\begin{quotation}
Take any $U_i$ and $U_j$ in $\fU$;
let $N$ and $N'$ be the restriction to $Z := U_i \cap U_j$
of the objects associated to $U_i$ and $U_j$ by $\beta$.
Let $X$ and $Y$ be any open sets
such that $\ol{X} \cap \ol{Y} = \emptyset$,
let $W = U'_1 \cup \cdots \cup U'_{p-1}$ for any integer $p$,
and let $f \colon N \to N'$ be any $\cB$ arrow.
Then there exists a $\cB$-arrow $f' \colon N \to N'$
with the following properties.
\begin{enumerate}
\item
The restriction of $f'$ to $X \cap Z$ is in $\cA$.
\item
The restrictions of $f'$ and $f$ to $Y \cap Z$
coincide.
\item
If the restriction of $f$ to $W \cap Z$ lies in $\cA$,
then the restriction of $f'$ to $W \cap Z$ also lies in $\cA$.
\end{enumerate}
\end{quotation}

Then the inclusion map $i \colon \cA \to \cB$ induces an isomorphism
$$ i_* \colon \vH^1(\calT,\cA) \stackrel{\cong}{\to} \vH^1(\calT,\cB).$$
\end{Lemma}

\begin{proof}
It is enough to show that
\begin{equation} \labell{map of Hone}
 i_* \colon \vH^1(\fU,\cA) \to \vH^1(\fU,\cB)
\end{equation}
is an isomorphism for every countable cover $\fU = \{U_i\}_{i=1}^\infty$
such that each $U_i$ intersects only a finite number of $U_j$'s.

\medskip \noindent
\emph{Proof that the map \eqref{map of Hone} is onto.} \\
Let $\beta \in Z^1(\fU,\cB)$ be a one cocycle.  We want to find
a zero cochain, $b \in C^0(\fU,\cB)$, such that
$b_V \circ \beta_{VU} \circ b_U\inv$ is in $\cA$ for all $U, V \in \fU$.
Let $\{ U_i' \}$ be a cover as in the statement of the lemma.

By induction, it  suffices  to prove
that for any one cocycle $\beta \in Z^1(\fU, \cB)$
such that the restriction of $\beta_{VU}$ to
$(U'_1 \cup \cdots \cup U'_{p-1}) \cap (U \cap V)$
is in $\cA$ for all $U, V \in \fU$,
we can find a zero cochain $b \in C^0(\fU,\cB)$
such that the restriction of $b_V \circ \beta_{VU} \circ {b_U}\inv$
to $(U'_1 \cup \cdots \cup U'_p) \cap (U \cap V)$
is in $\cA$ for all $U,V \in \fU$,  and
such that   $b_V$ is the identity for all $V \in \fU$ such that
$U_p \cap V = \emptyset$.
Since each $U_i$ intersects only a finite number of $U_j$'s,
the last condition ensures that for each  pair $U, V \in \fU$
the arrow $\beta_{VU}$ stabilizes after a finite number of steps.

Let $\beta \in Z^1(\fU,\cB)$ be a one cocycle such that
the restriction of $\beta_{VU}$ to
$(U'_1 \cup \cdots \cup U'_{p-1}) \cap (U \cap V)$
is in $\cA$ for all $U, V \in \fU$.
Let $Y$ be an open set such that
$Y \cup U_p = \calT$ and $\ol{Y} \cap \ol{U_p'} = \emptyset$.
Let $N_U$ denote the object associated by $\beta$ to $U$,
for each $U \in \fU$.
By the assumption of the lemma, for every $V \in \fU$
there exists a $\cB$-arrow
$\beta'_{VU_p} \colon N_{U_p} | _{U_p \cap V} \to N_V|_{U_p \cap V}$
with the following properties:
\begin{enumerate}

\item
The restriction of $\beta'_{VU_p}$ to $U'_p \cap V$ is in $\cA$.

\item
The restrictions of $\beta'_{VU_p}$ and $\beta_{VU_p}$
to $Y \cap (U_p \cap V)$ coincide.

\item
The restriction of  $\beta'_{VU_p}$ to $(U'_1 \cup \cdots \cup U'_{p-1})
\cap (U_p \cap V)$ is in $\cA$.
\end{enumerate}

By item (2), we may define $b_V \colon N_V \to N_V$ by
$b_V = {\beta'}_{VU_p}\circ \beta_{VU_p}\inv$
on $U_p \cap V$ and by $b_V = \mbox{id}$ on $Y \cap V$.
As required, if $U_p \cap V = \emptyset$ then $b_V = \mbox{id}$.

Now we claim that the restriction of $b_V \circ \beta_{VU} \circ b_U\inv$
to $(U'_1 \cup \cdots \cup U'_p) \cap (U \cap V)$ is in $\cA$,
for any $U, V \in \fU$.

By the definition of $b$, the restriction of
$b_V \circ \beta_{VU} \circ b_U\inv$
to $U_p \cap (U \cap V)$ is given by
$({\beta'}_{VU_p} \circ \beta_{VU_p}\inv) \circ \beta_{VU} \circ (\beta_{UU_p}
\circ {\beta'_{UU_p}}  \inv) = \beta'_{VU_p} \circ  {\beta'_{UU_p}} \inv$.
The restrictions of $\beta'_{VU_p}$ and $\beta'_{UU_p}$
to $(U_1' \cup \cdots \cup U_p') \cap U_p \cap U \cap V$ both lie in $\cA$,
by items (1) and (3).
Therefore, the restrictions of $b_V \circ \beta_{VU} \circ b_U\inv$ to
$(U_1' \cup \cdots \cup U_p') \cap U_p \cap (U \cap V)$ is also in $\cA$.

The restriction of $b_V \circ \beta_{VU} \circ b_U\inv$ to $Y \cap (U \cap V)$
coincides with the restriction of  $\beta_{VU}$ itself.
By the induction hypothesis, its restriction to
$(U'_1 \cup  \cdots \cup U'_{p-1})  \cap Y \cap (U \cap V)$
is in $\cA$.

Thus, the restriction of $b_V \circ \beta_{VU} \circ b_U\inv$ to the set
$(U_1 \cup \cdots \cup U_p') \cap U_p \cap (U \cap V)$ and to the set
$(U'_1 \cup  \cdots \cup U'_{p-1})  \cap Y \cap (U \cap V)$
lie in $\cA$.
Because the union of these sets is
$(U'_1 \cup \cdots \cup U'_p) \cap (U \cap V)$, we are done.

\medskip \noindent
\emph{Proof that the map \eqref{map of Hone} is one to one.} \\
Let $\alpha \in Z^1(\fU,\cA)$ be a one cocycle, and
let $b \in C^0(\fU,\cB)$ be a zero cochain such that
$b_V \circ \alpha_{VU} \circ b_U\inv$ is in $\cA$
for all $U, V \in \fU$.
We want to find a  zero cochain
$a \in C^0(\fU,\cA)$ such that
$a_V \circ \alpha_{VU} \circ a_U\inv = b_V \circ \alpha_{VU} \circ b_U\inv$
for all $U,V \in \fU$.
Let $\{ U_i' \}$ be a cover associated to
$i(\alpha) \in \vZ^1(\fU,\cB)$ as in the statement of the lemma.

By induction, it suffices to prove
that if we are given a zero cocycle $b \in C^0(\fU,\cB)$
such that $b_V \circ \alpha_{VU} \circ b_U\inv$ is in $\cA$
for all $U, V \in \fU$ and
the restriction of $b_V$ to $(U_1' \cup \ldots \cup U_{p-1}') \cap V$
is in $\cA$ for all $V \in \fU$, then we can find a zero cochain
$b' \in C^0(\fU,\cB)$ such that
$b'_V \circ \alpha_{VU} \circ {b'}_U\inv
  = b_V \circ \alpha_{VU} \circ b_U\inv$
for all $U, V \in \fU$,
the restriction of $b'_V$
to $(U'_1 \cup \cdots \cup U'_p) \cap V$
is in $\cA$ for all $V \in \fU$,
and such that if $V \cap U_p = \emptyset$
then $b'_V = b_V$.
Again, let $Y$ be such that
$Y \cup U_p = \calT$ and $\ol{Y} \cap \ol{U_p'} = \emptyset$.

Let $N_U$ denote the object associated by $i(\alpha)$ to $U$,
for each $U \in \fU$.
By the assumptions of the lemma, with $U_i = U_j = U_p$,
there exists a $\cB$-arrow $b_{U_p}' \colon N_{U_p} \to N_{U_p}$
with the following properties:
\begin{enumerate}
\item
The restriction of $b_{U_p}'$ to $U'_p$ is in $\cA$.
\item
The restrictions of $b_{U_p}'$ and $b_{U_p}$ to $Y \cap U_p$ coincide.
\item
The restriction of $b_{U_p}'$
to $(U'_1 \cup \ldots \cup U_{p-1}') \cap U_p$ is in $\cA$.
\end{enumerate}

By item (2), for every $V \in \fU$ we may
define a $\cB$-arrow  $b'_V \colon N_V \to N_V$
by $b'_V = b_V$ over $Y \cap V$ and
$b_V' = b_V \circ \alpha_{V U_p } \circ {b_{U_p}} \inv \circ {b'_{U_p}}
\circ \alpha_{V U_p }\inv$
on $U_p \cap V$.
As required, if $V \cap U_p = \emptyset$, then $b'_V = b_V$.

Now we claim that
$b_V' \circ \alpha_{VU} \circ b_U'{}\inv =
b_V \circ \alpha_{VU} \circ b_U \inv$ over $U\cap V$
for every $U,V \in \fU$.
Over $U \cap V \cap Y$ we have $b_V' = b_V$ and $b_U' = b_U$,
hence, $b_V' \circ \alpha_{VU} \circ b_U'{}\inv =
b_V \circ \alpha_{VU} \circ b_U \inv$.
Over $U \cap V \cap U_p$ we have
$b_V' \circ \alpha_{VU} \circ b_U'{}\inv =
(b_V' \circ \alpha_{V U_p } \circ b_{U_p}'{}\inv) \circ
(b_U' \circ \alpha_{U U_p } \circ b_{U_p}'{}\inv)\inv$
which, by the definition of $b_U'$ and $b_V'$, is equal to
$(b_V \circ \alpha_{V U_p} \circ b_{U_p} \inv) \circ
 (b_U \circ \alpha_{U U_p} \circ b_{U_p} \inv)\inv
= b_V \circ \alpha_{VU} \circ b_U \inv$.
Since $U_p \cup Y = \calT$, this proves that claim.

Finally, we claim that $b_V'$ is in $\cA$ over
$(U_1' \cup \ldots \cup U_p') \cap V$, for every $V \in \fU$.
Over $(U_1' \cup \cdots \cup U_p') \cap V \cap U_p$ we have
$b_V' = (b_V \circ \alpha_{V U_p } \circ {b_{U_p}} \inv)
\circ {b'_{U_p}} \circ \alpha_{V U_p }\inv$,
the first and third factors are in $\cA$ by assumption,
and the second is in $\cA$ by items (1) and (3).
Hence, their product is in $\cA$.
Over $(U_1' \cup \cdots \cup U'_p) \cap V \cap Y$,
$b'_V$ is equal to $b_V$.  Since $U'_p \cap Y = \emptyset$,
this set is the same as $(U_1' \cup \cdots \cup U_{p-1}') \cap V \cap Y$,
on which $b_V$ is in $\cA$ by the induction hypothesis.
Since $U_p \cup Y = \calT$, we are done.
\end{proof}

\section{Rigidification of $\Phi$-diffeomorphisms}
\labell{sec:rigidify Q to RQ}

In this section we prove that the sheaf $\hcQ$
of $\Phi$-diffeomorphisms
has the same first cohomology as the subsheaf $\cRQ$
of locally rigid $\Phi$-homeomorphisms.

\subsection{Definitions}

Consider two complexity one Hamiltonian $T$-manifolds $M$ and $M'$.
(As a special case we allow complexity one models.)
Let $\W \subseteq [0,1] \times M/T$ be an open subset,
and let
$$ G \colon \W \to M'/T $$
be an isotopy of $\Phi$-homeomorphisms.
(See Definition \ref{def:isotopy Phi homeo}.)
Let $g_t(\cdot) = G(t,\cdot)$.

\begin{Definition} \labell{def:isotopy Phi diffeo}
$G$ is an \textbf{isotopy of $\Phi$-diffeomorphisms}
if near each exceptional orbit $(t_0,E)$
it lifts to a smooth map $\tilde{G}$
from an open subset of $[0,1] \times M$ to $M'$
such that $\tilde{G}(t,\cdot)$ is a $\PhiT$-diffeomorphism
with its image for all $t$.
\end{Definition}

\begin{Definition} \labell{def:isotopy rigid}
$G$ is \textbf{rigid} at an exceptional orbit $(t_0,E)$
if for any pair of sub-grommets
$\psibar \colon D/T \to M/T$ and
$\psibar' \colon D/T \to M'/T$
whose images contain $E$ and $g_{t_0}(E)$
there exists a smooth function
$R \colon [0,1] \times \ft^* \to \Rbar$ such that
$$ {{\psibar}'}\inv \circ g_t \circ \psibar (y)
 = R(t,\Phibar_Y(y)) \cdot y $$
on some neighborhood of $(t_0,\psi\inv(E))$.
Note that both sub-grommets must have the same domain 
$D/T \subseteq Y/T$.
\end{Definition}

Note that a single $\Phi$-homeomorphism $g \colon W \to M'/T$,
where $W$ is an open subset of $M/T$, is locally rigid
if and only if it is rigid at every exceptional orbit in $W$.
(See Definition \ref{def:rigid on quotient} 
and Remark \ref{single Phi homeo}.)

\begin{Definition} \labell{rigidY}
Let $Y$ be a complexity one model.
Let $\W \subseteq [0,1] \times Y$ be an open subset.  A map
\begin{equation} \labell{G calW Y}
 G \colon \W \to Y 
\end{equation}
is an \textbf{isotopy} of $\PhiT$-diffeomorphisms
if it is smooth and if $g_t(\cdot) = G(t,\cdot)$
is a $\PhiT$-diffeomorphism (onto its image) 
for every $t \in [0,1]$.
(Here, we have the same model $Y$ in the domain and range.)
The isotopy $G$ is \textbf{rigid} at an exceptional orbit $(t_0,E)$
if it descends to an isotopy of $\Phi$-diffeomorphisms
$$ \ol{G} \colon \W/T \to Y/T $$
which is rigid at $(t_0,E)$. 
(See Definition \ref{def:isotopy rigid}.)
\end{Definition}

\begin{Lemma} \labell{rigid char}
An isotopy of $\PhiT$-diffeomorphisms \eqref{G calW Y}
is rigid at an exceptional orbit $(t_0,E)$ if and only if
there exists a smooth $T$-invariant function
$$ S \colon \W  \to R_Y $$
such that $G(t,y) = S(t,y) \cdot y$ 
on some neighborhood of $(t_0,E)$
and such that the composition 
of $S$ with the projection map $R_Y \to \Rbar$
is equal to the pullback via $\Phi_Y$ 
of a smooth function from $[0,1] \times \t^*$ to $\Rbar$.
\end{Lemma}

\begin{proof}
By definition, there exists a smooth function
$R \colon [0,1] \times \ft^* \to \Rbar$ so that 
$\overline{G}(t,y) = R(t,\ol{\Phi}_Y(y)) \cdot y$ 
on some neighborhood of $(t_0,E)$, where $\ol{G}\colon \W/T \to Y/T $ 
is the map induced by $G$.  
Let $\tilde{R} \colon [0,1] \times \t^* \to R_Y$
be a smooth map which is a lift of $R$ on some neighborhood
of $(t_0,\Phi_Y(E))$.
For each $(t,y)$ in a neighborhood of $(t_0,E)$,
the values $G(t,y)$ and $\tilde{R}(t,\Phi_Y(y)) \cdot y$  
are in the same $T$-orbit.
This implies that these maps differ on this neighborhood 
by a smooth $T$-invariant map from $\W$ to $T$,
by a theorem in \cite{HS} (see \cite[Theorem 4.12]{locun}).
We let $S$ be the product of this map with $\tilde{R} \circ \Phi_Y$.
\end{proof}

\subsection{Rigidification on $\mathbf{\C^{h+1}}$}

Let an $h$ dimensional group $H$ act on
$\C^{h+1}$ through an inclusion map $\rho = (\rho_0,\dots,\rho_h) 
\colon H \to (S^1)^{h+1}$.  
Let $L_\rho$ denote the group of $\R$-linear automorphisms
of $\C^{h+1}$ that preserve orientation, commute with the $H$-action,
and preserve the moment map $\Phi_H \colon \C^{h+1} \to \h^*$.
Let $R_\rho \subset L_\rho$ denote the subgroup of
unitary transformations that commute with the $H$-action, as in
Section~\ref{sec:def rigid}.

\begin{Lemma} \labell{rigid characterize}
$R_\rho$ is an $H$-equivariant smooth strong deformation retract 
of $L_\rho$.
\end{Lemma}

\begin{proof}
Let  $\eta_j = d\rho_j \in \h^*$ denote  the weights for the action.
Note that $\eta_j \neq 0$ if and only if $\rho_j(H)=S^1$.

Let us first assume that $\rho_0$ is equal to either $\rho_1$ or to
$\rho_1\inv$.
The complexity one assumption then implies that $\rho_j(H)=S^1$
and hence $\eta_j \neq 0$ for all $j$.  It also implies that $\rho_i$
is different from both $\rho_j$ and $\rho_j\inv$,
and hence they define non-isomorphic real representations of $H$
on $\C = \R^2$, for all $1 \leq i < j \leq h$.
The group of $\R$-linear transformations of $\C^{h+1}$
that commute with the $H$-action is, by Schur's lemma,
\begin{equation} \labell{A1 Ah}
   A_1 \times \ldots \times A_h,
\end{equation}
where $A_1$ consists of the $\R$-linear transformations of $\C^2$
which commute with the action of $H$ by $(\rho_0,\rho_1)$
and where $A_j = \Cc$ for $2 \leq j \leq h$.
We now have two sub-cases.
\begin{enumerate}
\item
Suppose that $\rho_0=\rho_1$. Then $A_1 = \GL(2,\C)$.
The moment map is
$$\Phi_H = \half \left( \eta_1 \left( |z_0|^2+|z_1|^2 \right)
  + \sum_{j=2}^h \eta_j |z_j|^2 \right).$$
The subgroup of \eqref{A1 Ah}
consisting of those elements that preserve (orientation and)
the moment map is $L_\rho = U(2) \times (S^1)^{h-1}$.
Thus, $L_\rho = R_\rho$.

\item
Suppose that $\rho_0 = \rho_1\inv$.
We apply the $\R$-linear transformation
$$(w_0,w_1,\ldots,w_h) = (\ol{z_0},z_1,\ldots,z_h).$$
In these new coordinates,
$H$ acts by $(\rho_1,\rho_1,\rho_2,\ldots,\rho_h)$,
and $A_1 = \GL(2,\C)$.
The moment map is
$$\Phi_H = \half \left( \eta_1 \left( - |w_0|^2 + |w_1|^2 \right)
   + \sum_{j=2}^h \eta_j |z_j|^2 \right).$$
The subgroup of \eqref{A1 Ah} consisting of those elements
that preserve (orientation and) the moment map is
$L_\rho = U(1,1) \times (S^1)^{h-1}$.
The group $R_\rho = (S^1)^{h+1}$ of rigid maps is a strong deformation
retract of $L_\rho$.
\end{enumerate}

Up to permutation, the only other case is that
$\rho_i$ is different from both $\rho_j$ and $\rho_j\inv$ for all $i \neq j$.
The group of real linear transformations of $\C^{h+1}$ that commute
with the $H$ action is, by Schur's lemma,
\begin{equation} \labell{A0 Ah}
  A_0 \times \ldots \times A_h,
\end{equation}
where each $A_i$ is the commutator of $\rho_i(H)$ in $\GL(2,\R)$.
The group of rigid maps is $R_\rho = (S^1)^{h+1}$.
We again distinguish between two sub-cases:
\begin{enumerate}
\item
Assume that  $\rho_0(H)$ is finite.
The complexity one condition then implies that $\rho_j(H) = S^1$,
so $A_j= \Cc$, for all $1 \leq j \leq h$.
We have $\eta_0 = 0$, and the moment map is
$$ \Phi_H = \half \sum_{j=1}^h \eta_j |z_j|^2 $$
with $\eta_j \neq 0$ for $1 \leq j \leq h$.

If $\rho_0(H) \not \subseteq \{1,-1\}$, then $A_1 = \Cc$.
If $\rho_0(H) \subseteq \{ 1, -1 \}$, then $A_1 = \GL(2,\R)$.
The subgroup of \eqref{A0 Ah} consisting of those elements that preserve
the orientation and the moment map is
$L_\rho = \Cc \times (S^1)^{h}$
or  $L_\rho = \GL^+(2,\R) \times (S^1)^{h}$, respectively.
In either case, $R_\rho = (S^1)^{h+1}$ is a strong deformation
retract of $L_\rho$.

\item
Up to permutation, the only other case is
where $\rho_i(H) = S^1$ for all $i$.
Then the group \eqref{A0 Ah} is $(\Cc)^{h+1}$.
The moment map is
$$ \Phi_H = \half \sum_{j=0}^h \eta_j |z_j|^2 $$
with $\eta_j \neq 0$ for all $j$.
The subgroup of \eqref{A0 Ah} consisting of
those elements of $(\Cc)^{h+1}$ that preserve
(the orientation and) the moment map is
$L_\rho = (S^1)^{h+1}$.  Thus $L_\rho = R_\rho$.
\end{enumerate}
\end{proof}

\subsection{Rigidification on a local model}

Consider a complexity one model
$$ Y = T \times_H \C^{h+1} \times \h^0 $$
with moment map $\Phi_Y \colon Y \to \t^*$.

\begin{Definition} \labell{Rplus action}
We let $\R_+ = \{ s \in \R \mid s \geq 0 \}$
act on
$Y = T \times_H \C^{h+1} \times \h^0$ by
$$ \mu_s \left( [\lambda,z,\nu] \right) 
   = [\lambda,sz,\nu], \quad s \in \R_+ . $$
\end{Definition}

\begin{Lemma} \labell{rigidify Y}
Fix a smooth function $s \colon [0,1] \to [0,1]$ such that
$s(t) = 0$ for $t$ in an open set containing $[\frac{1}{2},1]$
and $s(0) = 1$.
To a $T$-invariant open subset $W \subseteq Y$ we associate
an open subset
$\W \subseteq [0,1] \times Y$ by
$$\W := \{(t,w) \mid \mu_{s(t)}(w) \in W \} .$$
Let $g \colon W \to W' \subset Y$ be any $\PhiT$-diffeomorphism.
Then there exists an isotopy
of $\PhiT$-diffeomorphisms, $G \colon \W \to Y$,
with the following properties. Denote $g_t(\cdot) = G(t,\cdot)$.
\begin{enumerate}
\item $g_0 = g$.
\item $g_1$ is rigid at every exceptional orbit $E$.
\item If $g$ is rigid at $\mu_{s(t)}(E)$
      then $G$ is rigid at $(t,E)$,
      for any exceptional orbit $E$ and any $t$.
\end{enumerate}
(See Definitions \ref{rigidY} and \ref{Rplus action}.)
\end{Lemma}

\begin{proof}
Write $W = T \times_H V$, where $V \subseteq \C^{h+1} \times \h^0$
is open and $H$-invariant.
Since $g$ is a $\PhiT$-diffeomorphism, it locally has the form
\begin{equation} \labell{form}
 g ([ \lambda,z,\nu ]) = [\tau(z,\nu) \cdot \lambda, f(z,\nu),\nu],
\end{equation}
where $\tau \colon V \to T$ is smooth and $H$-invariant,
and where $f \colon V \to \C^{h+1}$ is smooth, and
for each $\nu \in \h^0$, the map $f(\cdot,\nu)$
is an $H$-equivariant diffeomorphism between open subsets of $\C^{h+1}$
that preserves the orientation and the moment map $\Phi_H$.
Because the origin is exceptional, $f(0,\nu) =0$ for all $\nu$.

Let $f_0(\nu) \colon \C^{h+1} \to \C^{h+1}$
be the $\R$-linear map obtained as the derivative of $f(\cdot,\nu)$
at the origin.  For $0 \leq t \leq \half$, define
\begin{equation} \labell{isotopy1}
  g_t ([\lambda,z,\nu]) = \begin{cases}
\left[ \tau(s(t) z,\nu) \cdot 
       \lambda, \frac{1}{s(t)} f(s(t)z,\nu), \nu \right] &
\text{ if } 0 < s(t) \leq 1 \\
\left[ \tau(0,\nu) \cdot \lambda, f_0(\nu)(z), \nu \right] &
\text{ if } s(t) = 0 .
\end{cases}
\end{equation}

The map $f_0(\nu)$ belongs to the group $L_\rho$ of $\R$-linear automorphisms
of $\C^{h+1}$ that preserve orientation, commute with the $H$-action,
and preserve the moment map $\Phi_H$.
Lemma \ref{rigid characterize} gives a smooth family of 
$H$-equivariant maps
$D_\sigma \colon L_\rho \to L_\rho$, for $0 \leq \sigma \leq 1$,
such that $D_0 = \id$, $\image D_1 = R_\rho$,
and $D_\sigma|_{R_\rho} = \id|_{R_\rho}$ for all $\sigma$,
where $R_\rho$ is the group of rigid maps of $\C^{h+1}$.

Let $\sigma \colon [\half,1] \to [0,1]$ be a smooth function
such that $\sigma(t) = 0$ for $t$ near $\half$ and $\sigma(1)=1$.
For $\half \leq t \leq 1$, define
\begin{equation} \labell{isotopy2}
g_t([\lambda,z,\nu])
= [ \tau(0,\nu) \cdot \lambda, D_{\sigma(t)}(f_0(\nu))(z), \nu].
\end{equation}

Note that whereas $\tau$ and $f$ in \eqref{form}
can be chosen in different ways, and can only be chosen
locally, $g_t$ only depends on $g$ and is therefore well defined.

Moreover, if $g$ is rigid at $\mu_{s(t_0)}(E)$, then 
there exists a $T$-invariant
smooth function $S\colon Y \to R_Y$ so that $g(y) = S(y) \cdot y$ on some
neighborhood of $\mu_{s(t_0)}(E)$, and so that the composition
of $S$ with the projection from $R_Y$ to $\Rbar$ is
the pull-back of a smooth function on $\ft^*$.
(See Lemma \ref{rigid char}.)
On a neighborhood of $(t_0,E)$,
$g_t$ is given by multiplication by $S(\mu_{s(t)}(y))$.
For $0 \leq t \leq \half$ this is a straightforward computation.
For $\half \leq t \leq 1$ this follows from the fact that
the deformation $D_{\sigma}$ fixes $R_\rho$.
Finally, by the formulas for $\mu_s$ and for $\Phi_Y$,
the function $(t,y) \mapsto S(\mu_{s(t)}(y))$
satisfies the properties in Lemma \ref{rigid char}.
\end{proof}

\subsection{Rigidification locally on $\mathbf{M/T}$}
\labell{subsec:rigidify M}

\begin{Lemma} \labell{Phi rigidify loc on M}
Let $(M,\omega,\Phi,\calT)$ and $(M',\omega',\Phi',\calT)$
be grommeted complexity one proper Hamiltonian $T$-manifolds.
Let $\psi \colon D/T \to M/T$ and $\psi' \colon D/T \to M'/T$
be sub-grommets with the same domain $D \subseteq Y$.
Suppose that $D/T$ is contractible.
Let $C \subset Y$ be a set of exceptional orbits
whose closure in $Y$ is contained in $D$.
Let $X \subseteq D$ be a set of exceptional orbits
such that $\mu_s(X) \subseteq X$ for all $0 \leq s \leq 1$.
(See Definition \ref{Rplus action}.)

For any $\Phi$-diffeomorphism $f \colon M/T \to M'/T$
that sends $\psi(C/T)$ to $\psi'(C/T)$ there exists
an isotopy of $\Phi$-diffeomorphisms
$F \colon [0,1] \times M/T \to M'/T$ with the following properties.
Denote $f_t = F(t, \cdot)$.
\begin{enumerate}
\item $f_0 = f$.
\item $f_1$ is rigid at every exceptional orbit in $\psi(C/T)$.
\item
$f_t = f$ outside $\psi(D/T) \cap f\inv \left( \psi'(D/T) \right)$.
\item
If $f$ is rigid at every orbit in $\psi(X/T)$,
then $F$ is rigid at every orbit in $[0,1] \times \psi(X/T)$.
\end{enumerate}
(See Definitions \ref{def:isotopy Phi diffeo} and \ref{def:isotopy rigid}.)
\end{Lemma}

\begin{Remark}
In fact, the isotopy $F$ depends ``smoothly"
on the function $f$.  This feature is relevant
for the study of the space of automorphisms
of a complexity one space. A similar situation occurs
in Section \ref{sec:rigidify P to RP} 
but not in Section \ref{sec:rigidify E}.
\end{Remark}

\begin{proof}[Proof of Lemma \ref{Phi rigidify loc on M}]
Once and for all, we choose a $T$-invariant smooth function
$$\rho \colon Y \to \R$$
such that
$$ \operatorname{support}(\rho) \subset D \quad \text{ and } \quad
\rho|_V \equiv 1$$
for some open neighborhood $V$ of $C$,
such that $\rho$ is a pullback of a smooth function on $\t^*$
on some neighborhood of the exceptional orbits.

Define open subsets of $Y$ by
\begin{align*}
 D_f & :=  \{ b \in D \mid f(\psi(b)) \in  \psi'(D/T) \} ,
\quad \mbox{and} \\
 D'_f & :=  \{ b' \in D \mid \psi'(b') \in  f(\psi(D/T)) \}.
\end{align*}

Because $D/T$ is contractible, $\psi\inv \circ f \circ \psi$
lifts to a $\PhiT$-diffeomorphism $g \colon D \to Y$.
In fact, it is enough to assume that $H^2(D/T,\Z) = 0$.
This is explained in our proof of Lemma 4.11 in \cite{locun},
following techniques of \cite{HS} and \cite{BM}.

Our maps fit into a commutative diagram:

$$ \begin{array}{ccccccc}
Y \supset & D & \supset D_f & \stackrel{g}{\to} &
                  D'_f \subset & D &  \subset Y \\
 & {\ss \psibar} \downarrow \phantom{\ss \psibar} & & & &
 \phantom{\ss \ol{\psi'}} \downarrow {\ss \ol{\psi'}} & \\
 & M/T & & \stackrel{f}{\to} & & M'/T . &
\end{array}$$
Moreover, $C \subset D_f \cap D'_f$
and $\ol{g}|_{C/T} = \id|_{C/T}$.

We now apply Lemma \ref{rigidify Y} with $W = D_f$
to obtain an isotopy $G(t,\cdot) = g_t(\cdot)$,
$$ g_t \colon \mu_{s(t)}\inv D_f \to \mu_{s(t)}\inv D_f'. $$
Let $\xi_t$ be the vector field which generates this isotopy.
This means that $\xi_t$ is a vector field defined on $\mu_{s(t)}\inv D_f'$
for each $t$, so that
$$ \frac{dg_t}{dt} (u) = \xi_t \circ g_t (u) $$
for each $u \in \mu_{s(t)}\inv D_f$ (as an equality in $T_{g_t(u)} Y$).
Let $\xi_t^\cutoff$  be the vector field
on $D'_f \cap \mu_{s(t)} \inv D_f'$ given by
\begin{equation} \labell{Yform again}
\xi_t^\cutoff(u)  := \rho(u) \, \rho( g\inv (u)) \,
\rho( \mu_{s(t)}(u)) \, \rho( g\inv (\mu_{s(t)}(u))) \, \xi_t(u).
\end{equation}

First, note that the support of $\xi_t^\cutoff$ in $Y$
is contained in the open set $D'_f \cap \mu_{s(t)} \inv D_f'$,
so $\xi_t^\cutoff$ extends to a smooth vector field $\xi_t^\cutoff$
on all of $Y$.
Construct an isotopy $h_t \colon D_f \to Y$ by solving the ordinary
differential equation
$$ \frac{dh_t}{dt} = \xi_t^\cutoff \circ h_t $$
with initial condition
$$ h_0 = g .$$
Since our cut-off functions are constant on orbits,
and since $g_t$  is $T$ equivariant, $h_t$ is also $T$ equivariant.
Similarly, $h_t$ respects the moment maps.

Second, note that there exists a closed subset of $D'_f$ such that
the support of $\xi_t^\cutoff$ is contained in this set for all $t$.
Therefore, each $h_t$ is a diffeomorphism from $D_f$ onto $D'_f$,
which coincides with $g$ on a neighborhood of the boundary
of $D_f$ in $Y$.
Using the sub-grommets $\psi$ and $\psi'$, the isotopy $h_t$
can be plugged back into $M/T$ to give an isotopy
of $\Phi$-diffeomorphisms $f_t \colon M/T \to M'/T$ such that
$f_t = f$ outside the image of $D_f$ and such that
${\psibar'}\inv \circ f_t \circ \psibar = \ol{h}_t$.

Third, note that $\xi_t^\cutoff$ coincides with $\xi_t$ on the set
$V'_f \cap \mu_{s(t)} \inv V'_f$, where
$V'_f = \{ v' \in V \mid \psi'(v') \in f(\psi(V)) \}$
(and $V$ is an open neighborhood of $C$ where $\rho=1$).
The intersection  of these sets for all $t \in [0,1]$
has a non-empty interior that contains $C$.
Hence, there exists an open neighborhood of $C$ in $D$
on which $h_t = g_t$ for all $t$.
Hence, $\ol{h}_1 = {\psibar'}\inv \circ f_1 \circ \psibar$
is rigid on $C/T$.

Finally, if $f$ is rigid at every orbit in $\psi(X/T)$,
then $g$ is rigid at every orbit in $X$.
Since $\mu_s(X) \subseteq X$ for all $0 \leq s \leq 1$, 
it follows from the third item of Lemma \ref{rigidify Y} 
that $G$ is rigid at every orbit in $[0,1] \times X$.
{}From Lemma \ref{rigid char} it follows that the time dependent
vector field $\xi_t$ has the following property. On a neighborhood 
of every exceptional orbit, $\xi_t$ 
is induced by a smooth $T$-invariant map from $[0,1] \times Y$ 
to the Lie algebra of $R_Y$, whose projection to
the Lie algebra of $\Rbar$ is the pullback of a smooth function
$[0,1] \times \t^* \to \Rbar$.
This implies that $\xi_t^\cutoff$  has
the same property (because the cut-off functions
are pullbacks from $\t^*$ near exceptional orbits).
Again by Lemma \ref{rigid char},
it follows that $h_t$ is rigid at every orbit in $[0,1] \times X$.
\end{proof}

\subsection{Rigidification globally on $\mathbf{M/T}$}
\labell{subsec:rigidify Q globally}

Let $(M,\omega,\Phi,\calT)$ be a complexity one Hamiltonian $T$-manifold.
For each point $\alpha \in \calT$,
let $\t_\alpha$ be the subspace of $\t$ spanned by all the
infinitesimal stabilizers to points in $\Phi\inv(\alpha) \subset M$.
Define an affine space $A_\alpha = \alpha + \t_\alpha^0$.
Let $p_\alpha \colon \t^* \to A_\alpha$
be the orthogonal projection determined by the fixed metric on $\t^*$.

\begin{Definition} \labell{def:orthogonal}
An open subset $V$ of $\calT$ is \textbf{orthogonal to the skeleton}
if and only if for each $\alpha \in \t^*$
there exists a neighborhood on which $V$ coincides with
a set of the form $p_\alpha\inv(V')$
for some open subset $V' \subseteq A_\alpha$.
\end{Definition}

\begin{Remark} \labell{refine cover}
For every open cover $\{ U_i \}$ of $\calT$
there exists an open cover $\{ U_i' \}$ such that
$\ol{U_i'} \subset U_i$ and $U_i'$ is orthogonal
to the skeleton for all $i$.
\end{Remark}

\begin{Proposition} \labell{need to rigidify Q to RQ}
Let $(M,\omega,\Phi,\calT)$ and $(M',\omega',\Phi',\calT)$ be
tall grommeted complexity one proper Hamiltonian $T$-manifolds.  Let
$$f \colon M/T \to M'/T$$
be a $\Phi$-diffeomorphism.  Let $W \subseteq \calT$ be an open subset
which is orthogonal to the skeleton. 
There exists an isotopy of $\Phi$-diffeomorphisms
$F \colon [0,1] \times M/T \to M'/T$ with the following properties.
Denote $f_t = F(t,\cdot)$.
\begin{enumerate}
\item $f_0 = f$.
\item $f_1$ is locally rigid.
\item If $f$ is rigid at every exceptional orbit 
      in $\Phi\inv(W)/T$, then $F$ is rigid 
      at every exceptional orbit in $[0,1] \times \Phi\inv(W)/T$.
\end{enumerate}
(See Definitions \ref{def:orthogonal}
\ref{def:isotopy Phi diffeo} and \ref{def:isotopy rigid}.)
\end{Proposition}

\begin{proof}
Define the \textbf{level} of an exceptional orbit
to be the dimension of its stabilizer.  Let $l$ be an integer.
Suppose that $f$ is rigid at all exceptional orbits
of level $> l$.  Let $W_l$ be the set of exceptional orbits
of level $l$ at which $f$ is rigid.
Let $K_l$ be the set of exceptional orbits of level $l$
at which $f$ is not rigid.  Note that $K_l$ is closed.

For each orbit $E$ in $K_l$ there exist sub-grommets
$\psibar_E \colon D_E/T \to M/T$ and 
$\psibar_E' \colon D_E/T \to M'/T$,
such that $\psibar_E([1,0,0]) = E$,
the domain $D_E/T$ is contractible,
and $f$ sends the exceptional orbits in $\psibar_E(D_E/T)$
to the exceptional orbits in $\psibar_E'(D_E/T)$.
Moreover, because $W$ is orthogonal to the skeleton,
we can choose the domains such that the set $X_E$
of exceptional orbits in $D_E \cap \Phi_E\inv(W)$ 
satisfies $\mu_s(X_E) \subseteq X_E$ for all $0 \leq s \leq 1$.
(See Definition \ref{Rplus action}.)
Since $K_l$ is closed, we can choose a collection of such sub-grommets
$$ \psibar_i \colon D_i/T \to M/T \quad \text{ and } \quad
   \psibar_i' \colon D_i/T \to M'/T ,$$
with domains $D_i \subset Y_i$
for $i = 1, \ldots, N \leq \infty$,
and such that the images $\psibar_i(D_i/T)$ form a locally finite
cover of $K_l$.  For each $i$, choose a subset $C_i$ of $D_i/T$
consisting of exceptional orbits of level $l$
whose closure in $Y_i/T$ is still contained in $D_i/T$,
and such that the images $\psibar_i(C_i)$ still cover $K_l$.

By induction on $i$, we construct a sequence of $\Phi$-diffeomorphisms
$$ f_i \colon M/T \to M'/T $$
such that $f_i$ is rigid at all points of 
$$ \Phi\inv(W) \cup W_l \cup \bigcup_{j \leq i} \ol{\psi}_j(C_j) $$
and a sequence of isotopies
$$ F_i \colon [0,1] \times M/T \to M'/T $$
such that $F_{i-1}(0,\cdot) = F_i(1,\cdot) = f_i(\cdot)$ for all $i$.

For $i=0$, set $f_0 = f$.  Given $f_{i-1}$, let
$$ X_i = \Phi_i\inv(W) \cup \psi_i\inv(W_l \cup \bigcup_{j<i} C_j) .$$
Notice that any subset of $D_i$ that consists of orbits of level $l$
is fixed by the $\R_+$ action $\mu_s$ on $D_i$ 
(see Definition \ref{Rplus action}).
It follows that $\mu_s(X_i) \subseteq X_i$ for all $0 \leq s \leq 1$.
Apply Lemma \ref{Phi rigidify loc on M} with the sub-grommet $\psibar_i$
and the subsets $C_i$ and $X_i$ to get an isotopy $F_i$, and set
$f_i(\cdot) := F(1,\cdot)$.

Each point $x \in M/T$ has a neighborhood $U$ and an $n=n(x)$
such that $F_j(s,y) =: f_\infty(y)$ is independent of $j$ and $s$
for all $j \geq n$ and $y \in U$.
Therefore, there exists an isotopy of $\Phi$-diffeomorphisms
$$ F \colon [0,1] \times M/T \to M/T $$
such that $F(0,\cdot) = f(\cdot)$,
$F(1 - \frac{1}{i+1}) = f_i(\cdot)$ for all $i$,
and $F(1,\cdot) = f_\infty(\cdot)$.
Note that $F(1,\cdot)$ is rigid at all exceptional orbits
of level $\geq l$. The result now follows by induction on $l$.
\end{proof}

We are ready to prove Proposition \ref{lemma RQ is hQ}.
We recall its statement:

\smallskip\noindent\textbf{Proposition \ref{lemma RQ is hQ}.}
{\it
The inclusion $\cRQ \subset \hcQ$ induces an isomorphism
$$ \vH^1(\calT,\cRQ) \cong \vH^1(\calT,\hcQ).$$
}\smallskip

\begin{proof}
The proposition follows from Lemma \ref{criterion for tech lemma}
once we show that the assumptions of this lemma are satisfied.

Let $\fU = \{ U_i \}$ be any cover
and $\beta \in \vZ^1 (\fU,\hcQ)$ be any one cocycle.
We can choose an open cover $\{ U_i' \}$
such that $\ol{U_i'} \subset U_i$ for each $i$
and such that $U_i' \cap U$ is an orthogonal set
with respect to $M_U$ for each $U \in \fU$,
where $M_U$ is the complexity one manifold associated to $U$ by $\beta$.

Take any $U_i$ and $U_j$ in $\fU$.  Let $M$ and $M'$ be the
restriction to $Z = U_i \cap U_j$ of $M_{U_i}$ and $M_{U_j}$.
Take any open sets $X$ and $Y$
such that $\ol{X} \cap \ol{Y} = \emptyset$, and let
$W = U_1' \cup \ldots \cup U_{p-1}'$, for any integer $p$.
Let $f \colon M/T \to M'/T$ be a $\Phi$-diffeomorphism.

Let $\rho \colon Z \to [0,1]$ be a smooth function
which vanishes on $Y \cap Z$ and is equal to one on $X \cap Z$.

Let $f_t$ be the isotopy obtained from Proposition
\ref{need to rigidify Q to RQ}. Then $f'(x) = f_{\rho(\Phi(x))} (x)$
fulfills the requirements (1)--(3) of Lemma \ref{criterion for tech lemma}.
\end{proof}

\section{Rigidification of local stretch maps}
\labell{sec:rigidify E}

In this section we prove that the sheaf of local stretch maps has
the same first cohomology as the subsheafs of locally rigid 
$\Phi$-homeomorphisms and of locally sb-rigid $\Phi$-homeomorphisms.

\subsection{Definitions}

\begin{Definition} \labell{def:Z}
Let $Y = T \times_H \C^{h+1} \times \h^0$ 
be a tall complexity one model with moment map $\Phi_Y$.
Its \textbf{associated painted surface bundle model} 
is the polyhedral subset
$$ Z = (\image \Phi_Y) \times \C \subseteq
\t^* \times \C,$$
and, for each exceptional orbit $E$ in $Y$,
a label on the point $p = F(E)$ of $Z$,
consisting of the isotropy representation at $E$.
Such a point is said to be \textbf{painted} by its label;
the set of labeled points in $Z$ is the \textbf{paint}.
\end{Definition}

\begin{remark}
We think of $Z$ as a bundle over $(\image \Phi_Y)$ with fiber $\C$.
\end{remark}

Let $Y$ and $Y'$ be tall complexity one models,
$Z$ and $Z'$ their associated surface bundle models,
and $F \colon Y \to Z$ and $F' \colon Y' \to Z'$
the trivializing homeomorphisms.
(See Section \ref{sec:stretch}.)
Let $W \subseteq Y/T$ and $\tilde{W} \subseteq Z$
be open subsets such that $F(W) = \tilde{W}$.
Let $G \colon [0,1] \times W \to Y'/T$
be an isotopy of $\Phi$-homeomorphisms
and $\tilde{G} \colon [0,1] \times \tilde{W} \to Z'$
be such that $\tilde{g}_t = F' \circ g_t \circ F\inv$,
where $g_t(\cdot) = G(t,\cdot)$
and $\tilde{g}_t(\cdot) = \tilde{G}(t,\cdot)$.

\begin{Definition} \labell{def:isotopy stretch}
$G$ and $\tilde{G}$ are \textbf{isotopies of stretch maps}
if there exists a function 
$R \colon [0,1] \times \Phi(W) \to S^1$ 
and an $S^1$-invariant function
$ \lambda \colon [0,1] \times \tilde{W} \to \R_{>0}$
such that
$$ \tilde{g}_t(\alpha,z) 
 = (\alpha , R(t,\alpha) \lambda(t,\alpha,z) \cdot z ) $$
for $(\alpha,z) \in \tilde{W} \subseteq \ft^* \times \C$.
\end{Definition}

\begin{Definition} \labell{def:isotopy sb rigid}
$\tilde{G}$ is \textbf{sb-rigid}
at a painted point $(t_0,p)$ if there exists a smooth function
$R \colon [0,1] \times \ft^* \to S^1$ such that
$$ \tilde{g}_t(\beta,z) = (\beta,R(t,\beta) \cdot z)$$
on some neighborhood of $(t_0,p) \in [0,1] \times Z
= [0,1] \times (\image \Phi_Y) \times \C$.
\end{Definition}

Now, let $(M,\omega,\Phi,\calT)$ and $(M',\omega',\Phi',\calT)$ 
be grommeted tall complexity one Hamiltonian $T$-manifolds.  
Let 
$$ F \colon [0,1] \times M/T \to M'/T $$
be an isotopy of $\Phi$-homeomorphisms.
Let $f_t(\cdot) = F(t,\cdot)$.

\begin{Definition} \labell{def:isotopy local stretch}
$F$ is an \textbf{isotopy of local stretch maps}
if for every exceptional orbit $E \in M/T$
and the pair of surface bundle grommets
$\varphi \colon B \to M/T$ 
and $\varphi' \colon B' \to M'/T$ 
whose images contain $E$ and $f_{t_0}(E)$,
the composition ${\varphi'}\inv \circ f_t \circ \varphi$
is an isotopy of stretch maps on some neighborhood
of $[0,1] \times E$.
(See Definition \ref{def:isotopy stretch}.)
\end{Definition}

Note that Definition \ref{def:isotopy local stretch}
is consistent with Definition \ref{def:stretch map}
for an isotopy that is independent of the parameter $t$.

\begin{Definition} \labell{def:sb rigid isotopy M}
$F$ is \textbf{sb-rigid} at an exceptional orbit $(t_0,E)$
if, for the pair of surface bundle grommets
$\varphi \colon B \to M/T$ and $\varphi' \colon B' \to M'/T$ 
whose images contain $E$ and $f_{t_0}(E)$,
the composition ${\varphi'}\inv \circ f_t \circ \varphi$
is sb-rigid on some neighborhood of $(t_0,E)$.
(See Definition \ref{def:isotopy sb rigid}.)
\end{Definition}

Note that a single $\Phi$-homeomorphism
$f \colon W \to M'/T$,
where $W$ is an open subset of $M/T$,
is locally sb-rigid if and only if it is sb-rigid
at every exceptional orbit in $W$.
(See Definition \ref{def:sb rigid} 
and Remark \ref{single Phi homeo}.)

\subsection{Rigidification on a local model}

Let $Y$ and $Y'$ be tall complexity one models,
$Z$ and $Z'$ their associated surface bundle models,
and $F \colon Y \to Z$ and $F' \colon Y' \to Z'$
the trivializing homeomorphisms.
(See Section \ref{sec:stretch}.)

\begin{Lemma} \labell{rigidify E on Y}
Let $W$ be an open subset of $Y/T$.  For any stretch map
$g \colon W \to Y'/T$
there exists an isotopy of stretch maps
$G \colon [0,1] \times W \to Y'/T$
with the following properties. let $g_t(\cdot) = G(t,\cdot)$.
\begin{enumerate}
\item $g_0 = g$.
\item $g_1$ is locally rigid.
\item If $g$ is rigid at $E$,
      then $G$ is rigid at $(t,E)$
	  for any exceptional orbits $E$ and any $t \in [0,1]$.
\end{enumerate}
(See Definitions \ref{def:isotopy stretch} and \ref{def:isotopy rigid}.)
\end{Lemma}

\begin{proof}
First, we define an $S^1$-invariant continuous function
$\tilde{\lambda} \colon F(W) \to \R_{>0}$
such that for every 
pair of canonical inclusions $\Lambda \colon D_E \to Y$ and
$\Lambda' \colon D_E \to Y'$ with the same domain we have
\begin{equation} \labell{good lambda}
 F' \circ \ol{\Lambda}' \circ \ol{\Lambda}\inv \circ F\inv(\alpha,z)
= (\alpha, \tilde{\lambda}(\alpha,z) \cdot z)
\end{equation}
on some neighborhood of the exceptional orbits
in $F(W \cap \ol{\Lambda}(D_E/T))$.

By Lemma \ref{rigid is stretch a}, 
for every such pair of canonical inclusions
there exists an $S^1$ invariant function $\tilde{\lambda}$
satisfying \eqref{good lambda}.
These functions agree on the intersections of their domains,
by Lemma \ref{canonical inclusion}
and because compositions of canonical inclusions are
canonical inclusions.
Hence, we can define $\tilde{\lambda}$ on a neighborhood
of the exceptional orbits in $F(W)$. We extend it arbitrarily to 
all of $F(W)$.

By assumption there exist $\lambda \colon F(W) \to \R_{>0}$ 
and $R \colon \Phi(W) \to S^1$ such that
$$F' \circ g \circ F \inv(\alpha,z) =
(\alpha, R(\alpha) \lambda(\alpha,z) \cdot z).$$
Define $g_t$ by
\begin{equation} \labell{this is gt}
F' \circ g_t \circ F \inv(\alpha,z) =
(\alpha, R(\alpha) \lambda_t(\alpha,z) \cdot z),
\end{equation}
where $\lambda_t = (1-t) \lambda + t \tilde{\lambda}.$

The function $\lambda_t$ is smooth on the complement of the
exceptional orbits, because so are the functions $\lambda$ 
and $\tilde{\lambda}$.
The maps $z \mapsto \lambda_t(\alpha,z) \cdot z$
have positive derivative everywhere because so do the maps
$z \mapsto \lambda(\alpha,z) \cdot z$ 
and $z \mapsto \tilde{\lambda}(\alpha,z) \cdot z$.
It follows that $\{ g_t \}$ is an isotopy of $\Phi$-homeomorphisms,
and, by \eqref{this is gt}, an isotopy of stretch maps.
\end{proof}

\subsection{Rigidification on $\mathbf{M/T}$}

\begin{Proposition} \labell{rigidify E to RQ}
Let $(M,\omega,\Phi,\calT)$ and $(M',\omega',\Phi',\calT)$
be grommeted tall complexity one proper Hamiltonian $T$-manifolds.
For any local stretch map $f \colon M/T \to M'/T$
there exists an isotopy of local stretch maps
$F \colon [0,1] \times M/T \to M'/T$ with the following properties.
Let $f_t(\cdot) = F(t,\cdot)$.
\begin{enumerate}
\item $f_0 = f$.
\item $f_1$ is locally rigid.
\item If $f$ is rigid at an exceptional orbit $E$,
then $F$ is rigid at $(t,E)$ for all $t$.
\end{enumerate}
(See Definitions \ref{def:isotopy local stretch} 
and \ref{def:isotopy rigid}.)
\end{Proposition}

\begin{proof}
We may assume that for each grommet $\psi \colon D \to M$
there exists a unique grommet $\psi' \colon D' \to M'$
such that $f$ sends the exceptional orbits in $\psi(D)$
to the exceptional orbits in $\psi'(D')$.

Let $W \subseteq D/T$ be a neighborhood of the exceptional orbits
on which the composition $g := {\psi'}\inv \circ f \circ \psi$
is a stretch map.  We apply Lemma \ref{rigidify E on Y}
to obtain an isotopy of stretch maps
$$g_t \colon W \to Y'/T$$
such that $g_1$ is locally rigid.

Let $\xi_t$ be the vector field on the complement
of the exceptional orbits in $g_t(W)$ for each $t$, such that
$$\frac {d {g}_t}{d t} = \xi_t \circ {g}_t.$$
Note that there exists a neighborhood of the exceptional orbits in $D'$
which is contained in $g_t(W)$ for all $t$.
Choose a function
$$ \rho \colon (Y'/T)|_\calT \to [0,1] $$
whose support is contained in $g_t(W) \cap D'/T$ for all $t$,
such that $\rho \equiv 1$ on a neighborhood of the exceptional orbits,
and such that the restriction of $\Phi_Y$ to the support of $\rho$
is proper.

The vector field
$$ \xi^\cutoff_t = \rho(u) \cdot \xi_t(u). $$
extends to a smooth vector field on the complement of the
exceptional orbits in $(Y'/T)|_\calT$, supported in $D'/T$.
Our goal is to find maps
$$ h_t \colon W \to (Y'/T)|_\calT $$
which satisfy the
ordinary differential equation
\begin{equation} \labell{ODE}
 \frac{d h_t}{dt} = \xi^\cutoff_t \circ h_t
\end{equation}
with initial condition $h_0 = g$.

Let $U' \subset (Y'/T)|_\calT$ be an open neighborhood of the exceptional orbits
on which $\rho(\cdot) \equiv 1$.
Let $U \subset W$ be an open neighborhood of the exceptional orbits
such that $g_t(U) \subseteq U'$ for  all $t \in [0,1]$.
If $x \in U$ then $g_t(x)$ solves \eqref{ODE}.

Let $V \subset W$ be an open set whose closure
does not contain any exceptional orbits
and such that $V \cup U = W$.
Let $V' \subset (Y'/T)|_\calT$ be an open set such
that $\ol{V'}$ does not contain any exceptional orbits
and  such that $g_t(V) \subseteq V'$
for all $t \in [0,1]$.
Let $\eta_t = \tilde{\rho} \cdot \xi_t^\cutoff$
where $\tilde{\rho} \equiv 1$ on $\ol{V'}$
and $\tilde{\rho} \equiv 0$ on a neighborhood of the exceptional orbits.
Because $L_{\eta_t} \Phi_{Y'} = 0$
and the restriction of $\Phi_{Y'}$ to the support of $\eta_t$
in $(Y' \ssminus Y'_\exc)|_\calT$ is proper,
there exists a solution $\tilde{h}_t$ to the differential equation
$ \frac{d \tilde{h}_t}{dt} = \eta_t \circ \tilde{h}_t$
with initial condition $\tilde{h}_0 = g$, defined on the complement
of the exceptional orbits in $W$.

Because $\eta_t = \xi_t^\cutoff$ on $V'$
and the restriction of $g_t$ to $V \cap U$ takes values in $V'$,
the restrictions to $V \cap U$ of $\tilde{h}_t$ and $g_t$ coincide.
Hence, we can define $h_t := g_t$ on $U$ and $h_t := \tilde{h}_t$
on $V$, and this solves \eqref{ODE}.

Because $h_t$ coincides with $g$
near the boundary of $W$ in $D/T$,
we can plug it back into the manifold $M/T$
to get an isotopy of local stretch maps with the required properties.
\end{proof}

We recall the statement of Proposition \ref{lemma RQ is E}.

\smallskip\noindent\textbf{Proposition \ref{lemma RQ is E}.}
{\it
The inclusion $\cRQ \subset \cE$ induces an isomorphism
$$ \vH^1(\calT,\cRQ) \cong \vH^1(\calT,\cE).$$
}\smallskip

Proposition \ref{lemma RQ is E} follows from
Proposition \ref{rigidify E to RQ}
in exactly the same way that Proposition \ref{lemma RQ is hQ}
followed from Proposition \ref{need to rigidify Q to RQ},
(except that the $U'_i$'s no longer need to be orthogonal
to the skeleton).

\subsection{sb-Rigidification}

\begin{Lemma} \labell{sb-rigidify E on Y}
Let $W$ be an open subset of $Z$.  For any stretch map
$$g \colon W \to  Z'$$
there exists an isotopy of stretch maps
$$ G \colon [0,1] \times W \to Z' $$
with the following properties. Let $g_t(\cdot) = G(t,\cdot)$.
\begin{enumerate}
\item $g_0 = g$.
\item $g_1$ is locally sb-rigid.
\item If $g$ is sb-rigid at $p$,
      then $G$ is sb-rigid at $(t,p)$,
	  for any painted point $p$ and any $t \in [0,1]$.
\end{enumerate}
(See Definitions \ref{def:isotopy stretch} 
and \ref{def:isotopy sb rigid}.)
\end{Lemma}

\begin{proof}
Let $\lambda \colon \t^* \times \C$ and $R \colon \t^* \to S^1$
be such that
$$ g (\alpha ,z) = (\alpha , R(\alpha) \lambda(\alpha,z) \cdot z).$$
Define 
$$ g_t (\alpha,z) = (\alpha , R(\alpha) \lambda_t(\alpha,z) \cdot z) $$
where $ \lambda_t = (1-t) \lambda + t \cdot 1$.
As in the proof of Lemma \ref{rigidify E on Y},
$G(t,\cdot) = g_t(\cdot) $ is an isotopy of stretch maps.
\end{proof}

\begin{Proposition} \labell{rigidify E to RP}
Let $(M,\omega,\Phi,\calT)$ and $(M',\omega',\Phi',\calT)$
be grommeted tall complexity one proper Hamiltonian $T$-manifolds.
For any local stretch map $f \colon M/T \to M'/T$
there exists an isotopy of local stretch maps
$F \colon [0,1] \times M/T \to M'/T$ with the following properties.
Let $f_t(\cdot) = F(t,\cdot)$.
\begin{enumerate}
\item $f_0 = f$.
\item $f_1$ is locally sb-rigid.
\item If $f$ is sb-rigid at an exceptional orbit $E$,
      then $\{ f_t \}$ is sb-rigid at $(t,E)$ for all $t$.
\end{enumerate}
(See Definitions \ref{def:isotopy local stretch} 
and \ref{def:sb rigid isotopy M}.)
\end{Proposition}

\begin{proof}
Proposition \ref{rigidify E to RP} follows from Lemma \ref{sb-rigidify E on Y}
in the same way that Proposition \ref{rigidify E to RQ} followed from
Lemma \ref{rigidify E on Y}.
\end{proof}

We recall the statement of Proposition \ref{lemma RP is E}:

\smallskip\noindent\textbf{Proposition \ref{lemma RP is E}.}
{\it
The inclusion $\cRP \to \cE$ induces an isomorphism
$$ \vH^1(\calT,\cRP) \cong \vH^1(\calT,\cE) .$$
}\smallskip

Proposition \ref{lemma RP is E} follows from 
Proposition \ref{rigidify E to RP}
in the same way that Proposition \ref{lemma RQ is E}
followed from Proposition \ref{rigidify E to RQ}.

\section{Rigidification of sb-diffeomorphisms}
\labell{sec:rigidify P to RP}

In this section we prove that the sheaf of sb-diffeomorphisms
has the same first cohomology as the subsheaf
of locally sb-rigid $\Phi$-homeomorphisms.

\subsection{Definitions}

\begin{Definition} \labell{def:family pi diffeo Z}
Let $Z$ be a painted surface bundle model
and $\W \subseteq [0,1] \times Z$
an open subset.  A map
$$ G \colon \W \to Z' $$
is an \textbf{isotopy of sb-diffeomorphisms}
if it sends each painted point to a painted point
with the same label, is smooth, 
and each $g_t (\cdot) = G(t , \cdot)$
is a diffeomorphism with its image.
\end{Definition}

\begin{Definition} \labell{def:sb isotopy on M}
Let $M$ and $M'$ be grommeted tall complexity one
Hamiltonian $T$-manifolds. An isotopy of $\Phi$-homeomorphisms
$$ F \colon [0,1] \times M/T \to  M'/T $$
is an \textbf{isotopy of sb-diffeomorphisms}
if for every pair of associated surface bundle grommets
$\varphi \colon B \to M/T$ and $\varphi' \colon B' \to M'/T$
the composition ${\varphi'}\inv \circ f_t \circ \varphi$
is an isotopy of sb-diffeomorphisms.
(See Definition \ref{def:family pi diffeo Z}.)
\end{Definition}

\subsection{Rigidification on a local model}

Let $Y$ and $Y'$ be tall complexity one models, and
let $Z$ and $Z'$ be their associated surface bundle models.

Let $Y = T \times_H \C^{h+1} \times \h^0$
and $\Phi_Y([t,z,\nu]) = \alpha + (\Phi_H(z) , \nu)$,
where the splitting $\t^* = \h^* \times \h^0$
is obtained from the metric on $\t$.
Then $\image \Phi_Y$ is the product of the affine space
$A_\alpha = \alpha + \h^0$ and the set $\image \Phi_H$,
so that
\begin{equation} \labell{splitting of Z}
 Z = A_\alpha \times \image \Phi_H \times \C.
\end{equation}

\begin{Definition} \labell{Rplus action on Z}
We define an $\R_+$-action on \eqref{splitting of Z} by
$$ \mu_s(q,\beta,z) = (q, s\beta, sz) \quad \text{ for } s \in \R_+. $$
\end{Definition}

\begin{Remark}
Definition \ref{Rplus action on Z} does not correspond to
Definition \ref{Rplus action} under the trivializing
homeomorphism $F \colon Y/T \to Z$.
\end{Remark}

\begin{Lemma}  \labell{rigidify Z}
Fix a smooth function  $s \colon [0,1] \to [0,1]$
such that $s(t) = 0$ for $t$ in an open set containing $[\half,1]$
and $s(0) = 1$.
To an open subset $W \subseteq Z$
we associate an open subset $\W \subseteq [0,1] \times Z$ by
$$ \W := \{ (t,w) \mid \mu_{s(t)}(w) \in W \} .$$
Let $g \colon W \to W' \subseteq Z$ be any sb-diffeomorphism.
Then there exists an isotopy of sb-diffeomorphisms,
$G \colon \W \to Z$, with the following properties.
Denote $g_t = G(t,\cdot)$.
\begin{enumerate}
\item $g_0 = g$.
\item $g_1$ is locally sb-rigid.
\item If $g$ is sb-rigid at $\mu_{s(t)}(p)$
      then $G$ is sb-rigid at $(t,p)$,
	  for any painted point $p$ and any $t$.
\end{enumerate}
(See Definitions \ref{def:family pi diffeo Z}
and \ref{def:isotopy sb rigid}.)
\end{Lemma}

\begin{proof}
Since $g$ is an sb-diffeomorphism, we have
$$ g (q,\beta,z) = (q,\beta,h(q,\beta,z)), $$
where $h \colon W \to \C$ is smooth,
and each map $h(q, \beta, \cdot)$ is an orientation preserving
diffeomorphism between open subsets of $\C$
that fixes the origin if $(q,\beta,0)$ is painted;
in particular, $h(q,0,0) = 0$.

Let
$$ h_0(q) \colon \h^* \times \C \to \C $$
be the $\R$-linear map obtained as the derivative of $h(q,\cdot,\cdot)$
at the origin.  The map $h_0(q)$ belongs to the group
of $\R$-linear maps from $\h^* \times \C$ to $\C$ of the form
$$B+A$$
where $A \colon \C \to \C$ is in $\GL_2^+(\R^2)$
and $B \colon \h^* \to \C$
is any linear map such that $B(\Phi_H(E)) = 0$ for each
exceptional orbit $E$ in $\C^{h+1}$.
This group, which we denote $L$,
strongly deformation retracts to the subgroup with $B=0$,
and, further, to the circle subgroup $R$ consisting of maps of the form
$$ (\beta,z) \mapsto \lambda z $$
for some $\lambda \in S^1$.  Choose a smooth family of maps
$D_\sigma \colon L \to L$, for  $0 \leq \sigma \leq 1$,
such that $D_0=\id$, $\image D_1 = R$,
and $D_\sigma|_{R} = \id |_{R}$ for all $\sigma$.

We first linearize. For $0 \leq t \leq \half$,
let $s=s(t)$, and define
$$ g_t (q,\beta,z) = \begin{cases}
(q,\beta,\frac{1}{s} h(q,s\beta,sz) ) & \text{ if } 0 < s \leq 1 \\
(q,\beta,h_0(q)(\beta,z) ) & \text{ if } s=0 .
\end{cases}$$

We now rigidify.
Let $\sigma \colon [\half,1] \to [0,1]$ be a smooth function
such that $\sigma(t)=0$ for $t$ near $\half$
and $\sigma(1)=1$.
For $\half \leq t \leq 1$, define
$$ g_t (q,\beta,z) = (q,\beta,D_{\sigma(t)} (h_0(q)) (\beta,z) ). $$
\end{proof}

\subsection{Rigidification locally on $\mathbf{M/T}$}

To rigidify, we  need to work with grommets in $M$ and $M'$
whose domains are the same.
The notion of ``sub-grommets" from Section~\ref{sec:def rigid}
is not good for this purpose, as the ``canonical inclusions" are not
sb-diffeomorphisms.
The ``surface bundle grommets" of Definition \ref{def:pi compatible}
are not good either, because their domains are prescribed.
We introduce the notion of ``surface bundle sub-grommets".

\begin{Definition}\labell{def:pi subgrommet}
Let $M$ be a grommeted tall complexity one proper Hamiltonian $T$-manifold.
Let $\varphi \colon B \to M/T$ be an associated surface bundle grommet
(Definition \ref{def:pi compatible}).
Let $E$ be a painted point in $B$
with associated surface bundle model $Z_E$.
Note that $B$ and $Z_E$ are both subsets of $\t^* \times \C$
which contain the point $E=(\pi(E),0)$.
Any sufficiently small neighborhood $B_E$ of $E$ in $Z_E$ 
is an open subset of $B$.
In this case, we call the restriction
$$ \varphi_E = \varphi|_{B_E} \colon B_E \to M/T $$
a \textbf{surface bundle sub-grommet} on $M$.
\end{Definition}

\begin{Lemma} \labell{rigidify locally in S}
Let $(M,\omega,\Phi,\calT)$ and $(M',\omega',\Phi',\calT)$ be
grommeted tall complexity one proper Hamiltonian $T$-manifolds.
Consider two associated surface bundle sub-grommets,
$\varphi \colon B \to M/T$ and $\varphi' \colon B \to M'/T$,
which have the same domain $B \subset Z = F(Y)$.
Let $C \subset Z|_\calT$ be a set of painted points
whose closure in $Z|_\calT$ is contained in $B$.
Let  $X \subset B$ be a set of painted points
such that $\mu_s(X) \subseteq X$ for all $0 \leq s \leq 1$.

For any sb-diffeomorphism $f \colon M/T \to M'/T$
that sends $\varphi(C)$ to $\varphi'(C)$, there exists an isotopy
of sb-diffeomorphisms
$F \colon [0,1] \times M/T \to M'/T$ with the following properties.
Denote $f_t = F(t,\cdot)$.
\begin{enumerate}
\item $f_0 = f$.
\item $f_1$ is locally sb-rigid.
\item $f_t(s) = f(s)$ for all $s$ outside
      $\varphi(B) \cap f \inv \varphi'(B)$.
\item If $f$ is sb-rigid at every exceptional orbit in $\varphi(X)$,
      then $F$ is sb-rigid at every exceptional orbit 
      in $[0,1] \times \varphi(X)$.
\end{enumerate}
(See Definition \ref{Rplus action on Z}, \ref{def:sb isotopy on M},
and \ref{def:sb rigid isotopy M}.)
\end{Lemma}

\begin{proof}
Once and for all, we choose a smooth function
$$\rho \colon Z \to \R$$
such that
$$ \operatorname{support}(\rho) \subset B \quad \text{ and } \quad
\rho|_V \equiv 1$$
for some open neighborhood $V$ of $C$,
such that 
$\rho$ is a pullback of a smooth function on $\t^*$
on some neighborhood of the painted points.

Define open subsets of $Z$ by
\begin{align*}
B_f & := \{ b \in B \mid f(\varphi(b)) \in  \varphi'(B) \} ,
\quad \text{and} \quad \\
B'_f & := \{ b' \in B \mid \varphi'(b') \in  f(\varphi(B)) \}.
\end{align*}
Our maps fit into a commutative diagram:
$$ \begin{array}{ccccccc}
Z \supset & B & \supset B_f & \stackrel{g}{\to} &
                  B'_f \subset & B &  \subset Z \\
& {\ss \varphi} \downarrow \phantom{\ss \varphi} & & & &
 \phantom{\ss \varphi'} \downarrow {\ss \varphi'} & \\
 & M/T & & \stackrel{f}{\to} & & M'/T &
\end{array}$$
Moreover, $C \subset B_f \cap B'_f$ and $g|_C = \text{identity}|_C$.

We now apply Lemma \ref{rigidify Z} with $W = B_f$
to obtain an isotopy
$$ g_t \colon \mu_{s(t)}\inv B_f \to \mu_{s(t)}\inv B_f'.$$
Let $\xi_t$ be the vector field which generates this isotopy.
This means that $\xi_t$ is a vector field defined on $\mu_{s(t)}\inv B_f'$
for each $t$, such that
$$ \frac{dg_t}{dt}(u) = \xi_t \circ g_t(u)$$
for each $u \in \mu_{s(t)}\inv B_f$ (as an equality in $T_{g_t(u)} Z$).
Let $\xi_t^\cutoff$  be the vector field
on $B'_f \cap \mu_{s(t)} \inv B_f'$ given by
\begin{equation} \label{Zform}
\xi_t^\cutoff(u)  := \rho(u) \, \rho( g\inv (u)) \,
\rho( \mu_{s(t)}(u)) \, \rho( g\inv (\mu_{s(t)}(u))) \, \xi_t(u).
\end{equation}

As in the proof of Lemma \ref{Phi rigidify loc on M},
we can construct an isotopy $h_t \colon B_f \to Z$ such that
$$ \frac{dh_t}{dt} = \xi_t^\cutoff \circ h_t $$
and
$$ h_0 = g ,$$
and plug this back into $M/T$. This gives an isotopy
of sb-diffeomorphisms $f_t \colon M/T \to M'/T$
which satisfies the conditions of the lemma. 
\end{proof}

\subsection{Rigidification globally on $\mathbf{M/T}$}

\begin{Proposition} \labell{need to rigidify P to RP}
Let $(M,\omega,\Phi,\calT)$ and $(M',\omega',\Phi',\calT)$
be grommeted tall complexity one proper Hamiltonian $T$-manifolds.  
Let 
$$f \colon M/T \to M'/T$$
be an sb-diffeomorphism.  Let $W \subseteq \calT$ be an open subset
which is orthogonal to the skeleton.
There exists an isotopy of sb-diffeomorphisms
$F \colon [0,1] \times M/T \to M'/T$ with the following properties.
Denote $f_t = F(t,\cdot)$.
\begin{enumerate}
\item $f_0 = f$.
\item $f_1$ is locally sb-rigid.
\item If $f$ is sb-rigid at each exceptional orbit in $W$,
      then $F$ is sb-rigid at each exceptional orbit
	  in $[0,1] \times W$.
\end{enumerate}
\end{Proposition}

\begin{proof}
This follows from Lemma \ref{rigidify locally in S}
in the same way that Proposition \ref{need to rigidify Q to RQ}
followed from Lemma \ref{Phi rigidify loc on M}.
\end{proof}

We recall the statement of Proposition \ref{lemma RP is hP}:

\smallskip\noindent\textbf{Proposition \ref{lemma RP is hP}.}
{\it
The inclusion $\cRP \subset \hcP$ induces an isomorphism
$$ \vH^1(\calT,\cRP) \cong \vH^1(\calT,\hcP).$$
}\smallskip

Proposition \ref{lemma RP is hP} follows from 
Proposition \ref{need to rigidify P to RP}
in the same way that Proposition \ref{lemma RQ is hQ}
followed from Proposition \ref{need to rigidify Q to RQ}.

\section*{Part V: Paintings}

\section{Surface bundles}
\labell{sec:S}

Let $(M,\omega,\Phi,\calT)$ be a tall complexity one space.
The map $\Phibar \colon M/T \to \calT$ is topologically
a fiber bundle over $\Phi(M)$ whose fibers are surfaces.
(See Proposition \ref{trivialize M mod T}.)
Because $\Phi(M)$ is a polyhedral subset of $\calT$,
the complexity one quotient $M/T$ is, topologically, a manifold with corners.
However, smoothly this is only true outside the exceptional orbits.
In this section we introduce \emph{painted surface bundles} over $\Phi(M)$.
Just like $M/T$, a painted surface bundle comes with a subset
that is ``painted" by isotropy data.
Unlike $M/T$, a painted surface bundle is a manifold with corners
\emph{everywhere}. 
It is relatively easy to determine whether two surface bundles are
isomorphic, and this will enable us to determine whether
two complexity one quotients are $\Phi$-diffeomorphic.

\begin{Definition}
A \textbf{skeleton} over an open set $\calT \subset \ft^*$ is a
topological space $S$ whose points are labeled by representations of
subgroups of $T$, together with a proper map $\pi \colon S \to \calT$.
This data must be locally modeled on the set of exceptional orbits
of a complexity one space in the following sense.
For each point $s \in S$, there exists a tall complexity
one model $Y = T \times_H \C^{h+1} \times \fh^0$ with
exceptional orbits $Y_\exc \subset Y/T$,
and a homeomorphism $\Psi$ from a neighborhood of $s$ to an open
subset of $Y_\exc$ which respects the labels and such
that $\ol\Phi_Y \circ \Psi = \pi$, where $\Phi_Y: Y \to \ft^*$
is the moment map.
\end{Definition}

\begin{Example}
The set $M_\exc$ of exceptional orbits in a
tall complexity one space is naturally a skeleton.
\end{Example}

Note that if $(S,\pi)$ is a skeleton over $\calT$,
then the map $\pi \colon S \to \calT$ is a local embedding.
Also, given an open subset $U \subset \calT$,
the restriction $S|_U := S \cap \pi\inv(U)$ is a skeleton over $U$.

If $(S,\pi)$ and $(S',\pi')$ are skeletons, an \textbf{isomorphism} 
from $S$ to $S'$  is a homeomorphism $i \colon S \to S'$
that sends each point to a point with the same isotropy data and
such that $\pi  = \pi' \circ i$.

Let $(S,\pi)$ be a skeleton. 
A function $\varphi \colon S \to \R$ is \textbf{smooth}
if for each point $s \in S$ there exists a neighborhood $U$
of $s$ in $S$ and a neighborhood $W$ of $\pi(s)$ in $\calT$
and a smooth function $\tilde{\varphi} \colon W \to \R$
such that $\tilde{\varphi} (\pi (x)) = \varphi(x)$ for all $x \in U$.
More generally, if $X$ is
a manifold (with corners), a map from $S$  to $X$ is \textbf{smooth}
if the pull-back of every smooth function on $X$ is a smooth function
on $S$.

\begin{Definition}
\labell{def:sb}
Let $\calT \subseteq \t^*$ be an open subset.
A \textbf{painted surface bundle} over $\calT$
is a manifold-with-corners $N$,\footnote
{\label{foot:corners}
A manifold with corners structure on $N$ is provided by an atlas
consisting of homeomorphisms $X_i \colon U_i \to \Omega_i$,
where $\{ U_i \}$ is an open cover of $N$ and $\Omega_i$
is an open subset of $\R_+^n$,
such that the transition maps $X_j \circ X_i \inv$ are smooth.
Here, a function on a subset of $\R^n$ is \textbf{smooth} if
it extends to a smooth function on an open subset of $\R^n$.
See \cite{corners}.
}
together with a proper map
$\pi \colon N \to \calT$ whose nonempty fibers are smooth oriented
surfaces, and a ``painted" subset $P \subset N$, whose points
are labeled by representations of subgroups of $T$,
subject to the following conditions.  First, every point in $N$ has a neighborhood
$U$ and a diffeomorphism $U \cong \pi(U) \times \mbox{(a disk)}$
which carries $\pi$ to the projection map to $\pi(U)$,
Second, $P$ is a skeleton, and the inclusion map
from $P$ to $N$ is smooth.
\end{Definition}

\begin{Definition}
\labell{sb:iso}
An \textbf{isomorphism} between painted surface bundles
is a diffeomorphism which respects the maps to $\t^*$,
the orientation on the fibers, and the paint.
\end{Definition}

\begin{Definition} \labell{assocS}
Let $(M,\omega,\Phi,\calT)$ be a grommeted tall complexity one
proper Hamiltonian $T$-manifold.
The \textbf{associated painted surface bundle} consists of
the following data:
\begin{enumerate}
\item
The topological manifold-with-corners $N_M = M/T$, together with
the map $\pi \colon N_M \to \calT$ that is induced by the moment map,
and the orientation on each fiber of $\pi$
obtained from the symplectic orientation of the reduced space
$\Phi\inv(\alpha)/T$.
\item
The manifold-with-corners structure on $N_M$
that is given by the following coordinate charts.
Choose arbitrary grommets whose images cover the complement
of the exceptional orbits in $M$ and are contained in this complement.
For each given grommet and each chosen grommet, take the associated
surface bundle grommet. (See Definition \ref{def:pi compatible}.)
\item
The subset $P$ of $N_M$ consisting of the exceptional orbits,
together with  a label for each $p \in P$
consisting of the isotropy representation
of the corresponding exceptional orbit in $M$.
We call this information the \textbf{paint}.
\end{enumerate}
\end{Definition}

The fact that the coordinate charts in item (2) give a well defined
smooth structure on $M/T$ follows from the facts that
the smooth structures given by the different grommets
coincide outside the set of exceptional orbits
(see Lemma \ref{F smooth}),
and that each exceptional orbit lies in the image of only one grommet.
The fact that we get a manifold with corners follows 
from \cite[Lemmas 4.7 and 7.1]{locun}.

Note that a map from $M/T$ to $M'/T$ is an sb-diffeomorphism
if and only if it is an isomorphism of the associated
painted surface bundles.
Here, $M$ and $M'$ are grommeted tall complexity one proper
Hamiltonian $T$-manifolds.

\begin{Definition} \labell{legal N}
A painted surface bundle $N$ over $\calT$ is \textbf{legal} if there
exists a cover $\{ W_\alpha \}$ of $\calT$ and for each $\alpha$
there exists a grommeted tall complexity one proper Hamiltonian $T$-manifold
$M_\alpha$ whose
associated painted surface bundle is isomorphic to $N|_{W_\alpha}$.
\end{Definition}

We introduce the sheaf $\cP$ of isomorphisms of surface bundles.
To each subset $U \subset \calT$
we associate a groupoid $\cP(U)$.
The objects in $\cP(U)$ are the legal painted surface bundles over $U$.
The arrows are the isomorphisms of painted surface bundles.

There is a natural inclusion map from $\hcP$ to $\cP$,
which associates to each grommeted complexity one space
its associated painted surface bundle and to each sb-diffeomorphism
the corresponding isomorphism of painted surface bundles.

\begin{Lemma}\labell{lemma hP is P}
The map $\hcP \to \cP$ induces an isomorphism in cohomology
$$\vH^1(\calT,\hcP) \cong \vH^1(\calT,\cP).$$
\end{Lemma}

\begin{proof}
This follows immediately from Lemma \ref{forgetful}.
\end{proof}

\section{Global objects: the exciting transition}
\labell{sec:exciting}

We are now ready, at last,
to leave the world of sheaves of groupoids, and return to
to the case where we have a single global object: a painted  surface bundle.

First, we extend Definition \ref{def:Phi homeo} as follows:

\begin{Definition} \labell{extend def:Phi homeo}
A $\mathbf{\Phi}$\textbf{-homeomorphism} 
between a complexity one quotient
and a painted surface bundle is a homeomorphism which 
respects the paint,
is a diffeomorphism off the paint, respects the maps to $\calT$,
and preserves orientations on the level sets of these maps.
\end{Definition}

\begin{Example}
Suppose that $M$ is grommeted.
The identity map from $M/T$ to the associated surface
bundle is a $\Phi$-homeomorphism.
\end{Example}

Let $\cQ$ be the sheaf of $\Phi$-diffeomorphisms as defined in
Section \ref{sec:grommets},
and let $\cP$ be the sheaf of isomorphisms of painted surface bundles
as defined in Section \ref{sec:S}.

By combining Lemmas \ref{lemma hQ is Q}, \ref{lemma RQ is hQ}.
\ref{lemma RQ is E}, \ref{lemma RP is E}, \ref{lemma RP is hP},
and \ref{lemma hP is P}, we get an isomorphism
\begin{equation} \labell{final iso}
\vH^1(\calT,\cQ) \cong \vH^1(\calT,\cP).
\end{equation}

By Lemma \ref{classifies over fU}, every legal painted surface bundle $N$
determines an element $[N]$ of $\vH^1(\calT,\cP)$.
Similarly, every tall complexity one proper Hamiltonian $T$-manifold $M$
determines an element $[M]$ of $\vH^1(\calT,\cQ)$.

\begin{Definition} \labell{def:N to M}
Let $(M,\omega,\Phi,\calT)$ be a tall complexity one
proper Hamiltonian $T$-manifold.
An \textbf{associated painted surface bundle}
is a legal painted surface bundle $N$ such that
$[M] \in \vH^1(\calT,\cP)$ and $[N] \in \vH^1(\calT,\cQ)$
correspond under the isomorphism \eqref{final iso}.
\end{Definition}

\begin{Proposition}\labell{main prop}
Let $(M,\omega,\Phi,\calT)$ be a tall complexity one
proper Hamiltonian $T$-manifold.  Up to isomorphism,
there exists a unique associated painted surface bundle $N$,
and it is $\Phi$-homeomorphic to $M/T$.

If $(M',\omega',\Phi,\calT)$ is another tall complexity one
proper Hamiltonian $T$-manifold,
and $N'$ is a painted surface bundle associated to $M'$,
then $M/T$ and $M'/T$ are $\Phi$-diffeomorphic
if and only if $N$ and $N'$ are isomorphic.
\end{Proposition}

\begin{proof}
Let $c_M \in \vH^1(\calT,\cP)$ be the class that corresponds
to $[M]$ under the isomorphism \eqref{final iso}.
Because the sheaf $\cP$ has gluable objects,
every element of $\vH^1(\calT,\cP)$ comes from a global object,
and this object is unique up to isomorphism.
(See Lemma \ref{classifies over fU}.)
In particular, there exists a legal painted surface bundle $N$
with $[N] = c_M$, and it is unique up to isomorphism.

By Lemma \ref{classifies over fU}, $M/T$ and $M'/T$ are $\Phi$-diffeomorphic
and only if $[M] = [M']$ in $\vH^1(\calT,\cQ)$.
Similarly, $N$ and $N'$ are isomorphic if and only if
$[N]=[N']$ in $\vH^1(\calT,\cP)$.
Because the map \eqref{final iso} that sends $[M]$ to $[N]$
and $[M']$ to $[N']$ is an isomorphism,
$[M]=[M']$ if and only if $[N]=[N']$.

It remains to prove that if $[M]$ corresponds to $[N]$
then $M/T$ is $\Phi$-homeomorphic to $N$.
Let $\fU$ be a countable cover of $\calT$
so that $M|_U$ can be grommeted for all $U \in \fU$, and so that
each open set in $\fU$ intersects only finitely many other open sets.
Let $a \in Z^1(\fU,\hcQ)$ be a cocycle such that for each $U \in \fU$ the
object $a_u$ is $M|_U$ with some choice of grommets, and the
arrows are the identity maps.  Clearly, under the
isomorphism $\vH^1(\calT,\hcQ) \to \vH^1(\calT,\cQ)$,
$[a]$ maps to $[M]$.

Let $[b] \in \vH^1(\fU,\hcP)$ be the image of $[a]$ under the isomorphism
of $\vH^1(\fU,\hcQ)$ with $\vH^1(\fU,\hcP)$.
Each of the isomorphisms composed to construct this isomorphism
was induced by an inclusion of two sheaves,
each of which is a  subsheaf of the sheaf $\hcH$ of $\Phi$-homeomorphisms
defined in Section \ref{sec:grommets}.
Therefore, $[a]$ and $[b]$ descend to the same class in $\vH^1(\fU,\hcH)$.
Thus, for each $U \in \fU$, there exists a $\Phi$-homeomorphism
from $M/T|_U$ to the object $b_U$, such that for each pair $U, V \in \fU$,
the associated arrow is $\beta_{UV} = f_U \circ f_V\inv$.

Under the inclusion $\vH^1(\fU,\hcP) \to \vH^1(\fU, \cP)$,
the class $[b]$ maps to $[N]$.  Hence, there exists $\Phi$-homeomorphisms
$g_u \colon b_U \to N|_U$ such that $g_U \circ g_V\inv = \id$ for each
pair $U, V \in \fU$. Composing with the $\Phi$-homeomorphisms
$M/T|_U \to b_U$, we get $\Phi$-homeomorphisms from $M/T|_U$ to $N|_U$
which coincide on overlaps.
\end{proof}

\section{Smooth Paintings}
\labell{sec:S embedding}

When $\calT$ is convex,
we can replace our painted surface bundle with a simpler object:
a smooth painting.

\begin{Definition} \labell{abstract painting}
Let  $(S,\pi)$ be a skeleton and
let $\Sigma$ be a closed oriented surface.
A \textbf{painting} is a  map $f \colon S \to \Sigma$
such that the map $(\pi,f) \colon S \to  \calT \times \Sigma$ is one to one.
Paintings $f \colon S \to \Sigma$ and $f'\colon S' \to \Sigma'$ are
\textbf{equivalent} if there exists an isomorphism
$i \colon S \to S'$ and a homeomorphism $\xi \colon \Sigma \to \Sigma'$
such that the compositions $\xi \circ f \colon S \to \Sigma'$
and $f' \circ i \colon S \to \Sigma'$ are homotopic through paintings.
\end{Definition}

One similarly defines smooth paintings and smooth equivalence of paintings.

\begin{Example} \labell{painting M}
Each tall complexity one space
naturally determines a painting $f \colon M_\exc \to \Sigma$
up to equivalence; see Section \ref{sec:statement and proof}.
\end{Example}

\begin{Proposition} \labell{bundle is painting}
Let $\calT \subseteq \t^*$ be an open set.
To every painted surface bundle $(N,P,\pi)$ over $\calT$ 
such that $\image \pi$ is convex and $\pi \colon N \to \image \pi$ is open,
one can associate a smooth equivalence class of paintings 
$f \colon P \to \Sigma$ such that
two such surface bundles $(N,P,\pi)$ and $(N',P',\pi')$
are isomorphic if and only if they have smoothly equivalent paintings
and $\image \pi = \image \pi'$.

Moreover, if $(M,\omega,\Phi,\calT)$ is a tall complexity one proper
Hamiltonian $T$-manifold and $(N,P,\pi)$ is an associated
painted surface bundle, then the associated paintings 
$M_\exc \to \Sigma$ and $P \to \Sigma$ are equivalent.
(See Section \ref{sec:statement and proof}.)
\end{Proposition}

The proof of the proposition relies on a ``trivialization"
of the surface bundles:

\begin{Lemma} \labell{trivialize N}
Let $\calT \subseteq \t^*$ be an open set.
Let $(N,P,\pi)$ be a painted surface bundle over $\calT$
such that $\image \pi$ is convex and $\pi \colon N \to \image \pi$
is open.  Then there exists a closed oriented surface $\Sigma$
and a smooth map $f \colon N \to \Sigma$ such that the map
$$ (\pi,f) \colon N \to (\image \pi) \times \Sigma $$
is a diffeomorphism.

Moreover, given any two such maps $f$ and $f'$
there exists a diffeomorphism
$\xi \colon \Sigma \to \Sigma$ such that $f$ is smoothly homotopic 
to $\xi \circ f'$ through maps which induce diffeomorphisms
from $N$ to $(\image \pi) \times \Sigma$.
\end{Lemma}

\begin{proof}
By the definition of a painted surface bundle, 
every point in $N$ has a neighborhood $U$
and a diffeomorphism $U \cong \pi(U) \times$ (a disc)
which carries $\pi$ to the projection map to $\pi(U)$.  
Because $\pi \colon N \to \image \pi$ is open, $\pi(U)$ is an
open subset of $\image \pi$.  
This implies that
every point in $\image \pi$ has a neighborhood $V$ in $\image \pi$ 
and a diffeomorphism $\pi\inv(V) \cong V \times \Sigma$ 
where $\Sigma$ is a closed oriented smooth surface 
which carries the map $\pi$ to the projection map
$V \times \Sigma \to V$.
The proof is similar to the proof of Ehresmann's lemma
(that a proper submersion is a fibration).
Because $\image \pi$ is contractible, the bundle $\pi \colon N \to \image \pi$
is trivial.
The second part of the lemma is proved exactly like 
Proposition \ref{trivialize M mod T}.
\end{proof}

\begin{proof}[Proof of Proposition \ref{bundle is painting}]
Let $f$ be as in Lemma \ref{trivialize N}.
By restricting  $f$ to the skeleton, we get 
a smooth painting, determined up to smooth equivalence.
Isomorphic painted surface bundles give isomorphic smooth paintings.

We need to show that if two painted surface bundles give rise
to smoothly equivalent paintings and have the same image in $\calT$, 
then they are isomorphic.

Given a skeleton $(S,\pi)$,
let $f_t \colon S \to \Sigma$ be a smooth homotopy through paintings.

Let $N = (\image \pi) \times \Sigma$,
and let $\pi$ be the natural projection to $\ft^*$.
Let $\hat{f}_t \colon S \to N$ be given by $\hat{f}_t(x) = (\pi(x),f_t(x))$.

Since $\hat{f}_t$ is smooth and one-to one,
for each $t$ and each $x \in S$
there exists a vector field $X_t$ on $N$, defined near $\hat{f}_t(x)$,
such that
\begin{equation} \labell{velocity}
\frac{d}{dt}\hat{f}_t(x) = X_t|_{\hat{f}_t(x)}.
\end{equation}
Note that $X_t$ is tangent to the fibers of $\pi$.
Using a partion of unity, one can obtain globally defined vector fields $X_t$
on $(\image \pi) \times \Sigma$ which are tangent to the fibers of $\pi$
and such that \eqref{velocity} holds.
We then integrate these vector fields to a family of diffeomorphisms
$g_t \colon (\image \pi) \times \Sigma \to (\image \pi) \times \Sigma$
such that $g_0 = \id$ and $\frac{d}{dt} g_t = X_t \circ g_t$.
Each $g_t$ preserves the fibers
of $\pi$, and $g_t(f(x)) = f_t(x)$ for all $x \in S$.
In particular, $g_1 \colon N \to N$ is an isomorphism
that respects the $\pi$ and such that $g_1(\hat{f_0}(x)) = f_1(x)$ for
all $x \in S$.

The last claim follows from the fact that, by Proposition \ref{main prop}, 
$M/T$ and $N$ are $\Phi$ homeomorphic.
\end{proof}

\section{Eliminating the smooth structure}
\labell{sec:smooth to continuous}
The final step is to show that, instead of working with smooth paintings,
we can simply work with (continuous) paintings.

\begin{Proposition} \labell{smooth to continuous}
Let $(S,\pi)$ be a skeleton, let $\Sigma$ be a closed oriented surface,
and consider paintings $f \colon S \to \Sigma$.
Every painting is equivalent to a smooth painting,
and if two smooth paintings are equivalent,
then they are smoothly equivalent.
\end{Proposition}

\begin{proof}
Embed $\Sigma$ into $\R^3$.
Choose an $\epsilon$ tubular neighborhood of $\Sigma$,
and let $r \colon U \to \Sigma$ be the associated normal retract.

We begin with the first claim.
Let $f \colon S \to \Sigma$ be a painting.

Choose a continuous function $\epsilon \colon S \to \R^+$ so that
$$\epsilon(s) < \epsilon \ \mbox{and} \
\epsilon(s) <  \frac{1}{4} ||f(s) - f(s')||
\ \forall \ s' \in S \ni \pi(s') = \pi(s) .$$

For all $s \in S$, choose a neighborhood $V_s \subset S$
such that
 $$||f(y) - f(s) ||  < \epsilon(y) \ \forall \ y \in V_s .$$
Let $\{U_\alpha\}$ be a locally finite refinement of $\{V_s\}$
with index assignment  $\alpha \mapsto s(\alpha)$.
Let $\lambda_\alpha$ be a smooth partition of unity
subordinate to $U_\alpha$.

Define $h \colon S \to \R^3$ by
$h(y) =  \sum_\alpha \lambda_\alpha(y)  f(s(\alpha));$
then $$|| h(y) -  f(y)||  \leq  \epsilon(y).$$

Since $\epsilon(y) < \epsilon$,  we can define $g_t \colon S \to \Sigma$
$$g_t(y) = r( (1-t) f(y) +   t \ h(y)) \ \forall t \in [0,1]$$
Clearly, $g_0 =  f$ and $g_1$ is smooth.
Moreover, $||g_t(y) - f(y)|| \leq 2\epsilon(y)$.
Therefore,  for any $y$ and $y'$ in $S$ such that $\pi(y) = \pi(y')$,
$$ || g_t(y) -  g_t(y') || \geq
 || f(y) - f(y') || -   ||g_t(y) - f(y) ||
-  || g_t(y') - f(y')|| > 0.$$
Therefore, $g_t$ is a painting.

We now prove the second claim.
Let $A$ denote $\{0,1\} \times S$.
Let $f \colon I \times S \to \Sigma$ be a homotopy of paintings such that
the restriction of $f$ to  $A$ is  smooth.

Choose a continuous function $\epsilon \colon I \times S \to \R^+$ so that
$$\epsilon(t,s) < \epsilon \ \mbox{and} \
\epsilon(t,s) <  \frac{1}{4} ||f(t,s) - f(t,s')||
\ \forall \ s' \in S \ni \pi(s') = \pi(s) .$$

For all $x \in I \times S$, choose a neighborhood $V_x \subset I \times S$
and a function $h_x \colon V_x \to \Sigma$ such that
\begin{itemize}
\item
If $x \in A$, then $h_x$ is a smooth local extension of $f|_{A \cap V_x}$.
\item
If $x \not \in A$, then $V_x \cap A = \emptyset$ and $h_x(y) = f(x)$.
\item For all $y \in V_x$, $|| f(y) - f(x) || < \epsilon(y)$ and
$||f(y) - h_x(y) ||  < \epsilon(y) .$
\end{itemize}

Let $\{U_\alpha\}$ be a locally finite refinement of $\{V_x\}$
with index assignment  $\alpha \mapsto x(\alpha)$.
Let $\lambda_\alpha$ be a smooth partition of unity
subordinate to $U_\alpha$.

Define $h \colon I \times S \to \R^3$ by
$h(y) =  \sum_\alpha \lambda_\alpha(y)  h_{x(\alpha)}(y);$
then
$$|| h(y) -  f(y)||  \leq  \epsilon(y).$$

Since $\epsilon(y) < \epsilon$,  we can define
$g \colon [0,1] \times S \to \Sigma$, by
$g(y) = r(h(y)).$

Clearly, $g|_A = f|_A$ and $g$ is smooth.
Moreover, $||g(y) - f(y)|| \leq 2\epsilon(y)$.
Therefore,  for any $y = (t,s)$ and $y' = (t,s')$ in $[0,1] \times S$
such that  $\pi(s) = \pi(s')$,
then
$$ || g(y) -  g(y') || \geq
 || f(y) - f(y') || -   ||g(y) - f(y) ||
-  || g(y') - f(y')|| > 0.$$
Therefore, $g$ is a homotopy through paintings.
\end{proof}

\section{Global existence up to $\Phi$-diffeomorphisms}
\labell{sec:globex}

In this paper, we have shown that certain invariants
determine a tall complexity one space up to isomorphism.
In our next paper, we will construct complexity one spaces
out of these invariants.
For future reference, 
we give here one step in this direction:  
given a skeleton $(S,\pi)$ and a closed convex subset 
$\Delta \subset \calT$, we show that if they locally come 
from complexity one spaces, then any painting 
$f \colon S \to \Sigma$ can be realized 
by gluing these spaces by $\Phi$-diffeomorphisms.

\begin{Proposition} \labell{technical}
Let $\calT$ be an open subset of $\ft^*$
and $\Delta \subset \calT$ a convex closed subset.
Let $f \colon S \to \Sigma$ be a painting, where $\Sigma$
is a closed oriented surface of genus $g$
and $(S,\pi)$ is a skeleton over $\calT$.
Let $\fU$ be an open cover of $\calT$.
For each $U \in \fU$, let $(M_U,\omega_U,\Phi_U,U)$
be a connected complexity one proper Hamiltonian $T$-manifold
of genus $g$ over $U$,
so that $\image \Phi_U = U \cap \Delta$
and so that the set of exceptional orbits $(M_U)_\exc$
is isomorphic to the restriction $S|_U := S \cap \pi\inv(U)$.

Then, after possibly refining the open cover,
one can associate to each $U$ and $V$ in $\fU$ a $\Phi$-diffeomorphisms
$h_{V U} \colon M_U/T|_{U \cap V} \to  M_V/T|_{U \cap V}$
such that $h_{W V} \circ h_{V U} = h_{W U}$,
and such that the following holds.

If $(M,\omega,\Phi,\calT)$ is a complexity one space  
such that for every $U \in {\mathcal U}$ there
exists a $\PhiT$-diffeomorphism $\lambda_U \colon M|_U \to M_U$
so that $h_{V U}$ is the map induced by
the composition $\lambda_V \circ (\lambda_U)\inv$,
then the painting associated to $M$ is equivalent to $f$.
\end{Proposition}

The proof of Proposition \ref{technical} will use the 
fact that, locally, a painted surface bundle is uniquely
determined by its skeleton, its image, and its genus:

\begin{Lemma}[Local uniqueness of surface bundles]
\labell{local uniqueness sb}
Let $(N,P,\pi)$ and $(N',P',\pi')$
be two painted surface bundles over an open neighborhood
of a point $\alpha$ in $\t^*$.
Suppose that they have the same genus,equivalent skeletons,
and the same image in $\t^*$.  Then there exists a neighborhood
$U$ of $\alpha$ such that the restrictions 
$(N|_U,P|_U,\pi|_U)$ and $(N'|_U,P'|_U,\pi'|_U)$ are isomorphic.
\end{Lemma}

\begin{proof}
Let $f \colon P \to \Sigma$ be a smooth painting
associated to $N$ as in Proposition \ref{bundle is painting}.
By the local normal form theorem,
there exists a neighborhood $U$ of $\alpha$ so
that $\pi$ identifies each component of $P|_U$ 
with a subset of $\t^*$ which is star shaped about $\alpha$.
Hence, $f|_U$ is smoothly homotopic to a trivial painting;
on each component, simply define $f_t(x) = f(t \alpha + (1 - t) x)$.
Applying a similar argument to $N'$, $N|_U$ and $N'|_U$
have smoothly equivalent paintings.
By Proposition \ref{bundle is painting}, $N|_U$ and $N'|_U$ 
are isomorphic.
\end{proof}

We can now prove our main proposition.

\begin{proof}[Proof of Proposition \ref{technical}]
By Proposition \ref{smooth to continuous},
after possibly passing to an equivalent painting,
we may assume that the painting $f \colon S \to \Sigma$ is smooth.
Let $N = \Delta \times \Sigma$, let $\pi \colon N \to \calT$
be the projection to the first coordinate, and let $P \subset N$
be the subset defined by $P = \{ (\pi(s),f(s)) \mid s \in S \}$,
with the labels inherited from $S$.  
For each $U \in {\mathcal U}$, since $M_U$ is tall there
are non-exceptional orbits in  each nonempty moment map fiber. 
Hence, the local normal form theorem, together with the stability
of the moment map image, imply 
that $\image \Phi_U$ is a manifold with corners. (See Lemma 4.7
in \cite{locun}.)
Since, by assumption, $\image \Phi_U = U \cap \Delta$,
this implies that $N$ is a manifold with corners.
By Definition \ref{def:sb}, $(N,P,\pi)$ is a painted surface bundle.

Let $\alpha \in \calT$ be any point.
Let $(N',P',\pi')$ be a painted surface bundle associated
to $(M_U,\omega_U,\Phi_U,U)$, as in Proposition \ref{main prop},
where $\alpha \in U \in {\mathcal U}.$  
By Lemma \ref{local uniqueness sb} there exists
a neighborhood $W$ of $\alpha$ such that
$(N|_W,P|_W,\pi|_W)$ is isomorphic to $(N'|_W,P'|_W,\pi'|_W)$,
so that $(N|_W,P|_W,\pi|_W)$ is a surface bundle
associated to $(M_U|_W,\omega_U|_W,\Phi_U|_W,W)$.
By Definition \ref{legal N}, the surface bundle $(N,P,\pi)$ is legal.

Hence, by Lemma \ref{classifies over fU}, it
gives rise to an element in the first cohomology $H^1(\calT,\cP)$,
where $\cP$ is the sheaf of isomorphisms of legal surface bundles
defined at the end of Section \eqref{sec:S}.
By the isomorphism \eqref{final iso}, this, in turn, comes
from a cohomology class in $H^1(\calT,\cQ)$.
Let $\{ h_{\alpha \beta} \}$
be a cocycle that represents this class.

Let $(M,\omega,\Phi,\calT)$ be as stated in the proposition.
By the definitions, the element of $H^1(\calT,\cP)$
associated to $M$ is exactly the one represented by
the cocycle $\{ h_{\alpha\beta} \}$.
By Definition \ref{def:N to M},
this means that $(N,P,\pi)$ is an associated surface bundle
to $(M,\omega,\Phi)$.  By Proposition \ref{main prop},
this implies that $M/T$ is $\Phi$-homeomorphic to $N$.
Therefore, the paintings associated to $M/T$ and $N$ are equivalent.
\end{proof}


\begin{thebibliography}{LMTW}

\bibitem[AH]{ah-ha} K. Ahara and A. Hattori,
\emph{ 4 dimensional symplectic $S^1$-manifolds admitting moment map},
J.\ Fac.\ Sci.\ Univ.\ Tokyo Sect.\ IA, Math.\ \textbf{38} (1991), 251--298.

\bibitem[Au1]{audin:paper} M. Audin,
\emph{Hamiltoniens p\'{e}riodiques sur les vari\'{e}t\'{e}s symplectiques
compactes de dimension 4}, G\'{e}om\'{e}trie symplectique et
m\'{e}canique, Proceedings 1988, C.\ Albert ed., Springer Lecture Notes in
Math.\ \textbf{1416} (1990).

\bibitem[Au2]{audin:book} M. Audin,
\emph{The topology of torus actions on symplectic manifolds},
Progress in Mathematics \textbf{93}, Birkhäuser Verlag, Basel, 1991.


\bibitem[BM]{BM} M.\ Boucetta and P.\ Molino,
\emph{G\'eom\'etrie globale des syst\`emes hamiltoniens compl\`ement
int\'egrables: fibrations lagrangiennes singuli\`eres et coordonn\'ees
action-angle \`a singularit\'es.}
(Frence. English summary)
\emph{[Global geometry of completely integrable Hamiltonian systems:
Lagrangian singular foliations and action-angle variables with 
singularities]}
C.\ R. Acad.\ Sci.\ Paris S\'er.\ I Math.\ \textbf{308} (1989),
no.\ 13, 421--424.

\bibitem[Br]{Br} J.-L. Brylinski,
\emph{Loop Spaces, Characteristic Classes, and Geometric Quantization},
Birkh\"auser Boston 1993.

\bibitem[C]{river} R. Chiang, \emph{Complexity one Hamiltonian
$\SU(2)$ and $\SO(3)$ actions}, Amer.\ J.\ Math., to appear.


\bibitem[De]{De} T. Delzant,
\emph{Hamiltoniens p\'{e}riodiques et image convexe de l'application
moment}, Bull.\ Soc.\ Math.\ France \textbf{116} (1988), 315--339.

\bibitem[DH]{corners}
A. Douady and L. H\'erault, appendix to:
A.\ Borel and J.-P.\ Serre, \emph{Corners and arithmetic groups},
Comment.\ Math.\ Helv.\ \textbf{48} (1973), 436--491.



\bibitem[GLS]{GLS} V.\ Guillemin, E.\ Lerman and S.\ Sternberg,
\emph{Symplectic Fibrations and Multiplicity Diagrams}, 
Cambridge Univ.\ Press, 1996.


\bibitem[GS2]{GS:normal} V. Guillemin and S. Sternberg,
\emph{A normal form for the moment map},
Differential geometric methods in mathematical physics,
(S. Sternberg, Ed.) Reidel, Dordrecht, Holland, 1984.


\bibitem[HS]{HS}
A.\ Haefliger  and  E.\ Salem, \emph{Actions of tori on orbifolds},
Ann.\ Global Anal.\ Geom.\ \textbf{9} (1991), 37--59.

\bibitem[Hu]{husemoller}
D. Husemoller, emph{Fiber Bundles}, Third Edition,
Springer-Verlag, New York, Inc.\ 1994.


\bibitem[K]{karshon:periodic}
Y.\ Karshon, \emph{Periodic Hamiltonian flows on four dimensional manifolds},
Mem.\ Amer.\ Math.\ Soc.\ \textbf{672} (1999).

\bibitem[KT]{locun}
Y.~Karshon and S.~Tolman,
\emph{Centered complexity one Hamiltonian torus actions},
Trans.\ Amer.\ Math.\ Soc.\ \textbf{353} (2001), 4831--4861.


\bibitem[LMTW]{LMTW} E. Lerman, E. Meinrenken, S. Tolman, and C. Woodward,
\emph{Nonabelian convexity by symplectic cuts},
Topology \textbf{37} (1998), 245--259.

\bibitem[L]{L} H. Li, \emph{Semi-free Hamiltonian circle actions on
6-dimensional symplectic manifolds}, Trans.\ Amer.\ Math.\ Soc.\
\textbf{355} (2003), no.~11, 4543--4568.

\bibitem[M]{marle}
C.\ M. Marle, { Mod\`{e}le d'action hamiltonienne d'un groupe de
Lie sur une vari\'{e}t\'{e} symplectique}, \emph{Rendiconti del
Seminario Matematico} \textbf{43} (1985), 227--251, Universit\`{a}e
Politechnico, Torino.


\bibitem[Mo]{moser}
J. Moser, \emph{On the volume elements on a manifold},
Trans.\ Amer.\ Math.\ Soc.\ \textbf{120} (1965), 286--294.

\bibitem[Sch]{schwarz:IHES} G.\ W. Schwarz,
\emph{Lifting smooth homotopies of orbit spaces},
I.H.E.S.\ Publ.\ Math.\ \textbf{51} (1980), 37--135.

\bibitem[T1]{T1}  D. A. Timash\"ev,
\emph{$G$-manifolds of complexity $1$},
(Russian) Uspekhi Mat.\ Nauk \textbf{51} (1996), no.\ 3 (309), 213--214;
English translation: Russian Math.\ Surveys \textbf{51} (1996), 567--568.

\bibitem[T2]{T2}  D. A. Timash\"ev,
\emph{Classification of $G$-manifolds of complexity $1$},
(Russian) Izv.\ Ross.\ Akad.\ Nauk Ser.\ Mat.\ \textbf{61} (1997), 127--162;
English translation: Izv.\ Math.\ \textbf{61} (1997), 363--397.

\bibitem[T]{T}
S. Tolman, \emph{Examples of non-K\"ahler Hamiltonian torus actions},
Invent.\ Math.\ \textbf{131} (1998), no.\ 2, 299--310.

\end{thebibliography}
\end{document}